\documentclass[10pt,a4paper]{amsart} 
\usepackage{amsfonts,amssymb,amsthm,amsmath,ifpdf,stmaryrd,amsrefs,dsfont,esint,mathrsfs,dsfont,xparse,arydshln}
\usepackage[hidelinks]{hyperref}
\usepackage{enumitem} 

\usepackage{tikz}
\usetikzlibrary{matrix,arrows}
\usepackage{tikz-cd}

\usepackage[overload]{textcase}

\ifpdf
  \usepackage[final,expansion=alltext,protrusion=alltext]{microtype} 
\fi

\usepackage{xargs} 
\usepackage{mathtools} 
\usepackage{ifthen} 

\theoremstyle{plain}

\newtheorem*{Th*}{Theorem}
\newtheorem*{ThA}{Theorem A}
\newtheorem*{ThB}{Theorem B}

\newtheorem*{Cor*}{Corollary}

\theoremstyle{definition}

\theoremstyle{remark}

\newif\ifTNS 
\TNStrue

\def\printtheoremname#1{\csname#1name\endcsname}
\def\printtheoremnames#1{\csname#1names\endcsname}

\def\thmref#1#2{\printtheoremname{#1}\ifTNS~\fi\ref{#1:#2}}

\def\uc#1#2{\MakeUppercase{#1}{#2}} 

\newcommand{\DefTheorem}[2]{\newenvironmentx{#1}[2][1=\empty,2=\empty]{%
    \ignorespaces%
    \ifx##2\empty%
      \begin{#2}%
    \else%
      \begin{#2}[{\uc##2}]%
    \fi%
    \ifx##1\empty%
      {}%
    \else%
      \label{#1:##1}%
    \fi%
    \ignorespaces}{\end{#2}\ignorespacesafterend}}

\newcommand{\prfof}[2]{\protect{Proof of~\thmref{#1}{#2}}}

\makeatother

\DefTheorem{Th}{theorem}
\DefTheorem{Prop}{proposition}
\DefTheorem{Cor}{corollary}
\DefTheorem{Lem}{lemma}
\DefTheorem{Def}{definition}
\DefTheorem{Rem}{remark}
\DefTheorem{Par}{para}
\DefTheorem{Not}{notation}
\DefTheorem{Exer}{exercise}
\DefTheorem{Ex}{example}
\DefTheorem{Cons}{construction}
\DefTheorem{Scho}{scholium}

\newenvironment{Par*}{\ignorespaces\noindent\ignorespaces}{\ignorespacesafterend}

\numberwithin{equation}{section}

\DeclareDocumentCommand\fibint{ O{\empty} m O{\empty} }{\!\sideset{_{#1}}{_{#2}^{#3}}{\fint}}
\DeclareDocumentCommand\tfibint{ O{\empty} m O{\empty} }{\!\sideset{_{#1}}{_{#2}^{#3}}{\textstyle\fint}}

\newcommand\Define[2][\empty]{\ignorespaces%
  \emph{#2}}%

\tikzset{
  commutative diagrams/.cd,
  arrow style=tikz,
  diagrams={>=stealth},
  shift up/.style={
    to path={([yshift=#1]\tikztostart.east) -- ([yshift=#1]\tikztotarget.west) \tikztonodes}},
  shift up left/.style={
    to path={([yshift=#1]\tikztostart.west) -- ([yshift=#1]\tikztotarget.east) \tikztonodes}},
  mathdouble/.style={-,double equal sign distance}
}

  {\ignorespaces\end{tikzpicture}%
  \end{minipage}\\%
  \ignorespacesafterend}

\def\ger{\mathfrak}

\newcommand\CategoryTypeface{\mathbf}
\def\cat{\CategoryTypeface}

\newcommand\SheafTypeface{\mathcal}
\def\sh{\SheafTypeface}

\def\DMO{\DeclareMathOperator}

\newcommand\ev{{\bar 0}}
\newcommand\odd{{\bar 1}}

\newcommand{\defi}{\coloneqq}     

\def\diff#1^#2{\ensuremath{\partial_{#1}^{#2}}}

\def\der#1/#2{\ifthenelse {\equal{#1}{}}
              {\ensuremath{\partial_{#2}{#1}}}
              {\ensuremath{\frac{\partial #1}{\partial #2}}}
        }

\def\derf#1/#2{\ifthenelse  {\equal{#1}{}}
              {\ensuremath{\frac{\partial #1}{\partial #2}}}
              {\ensuremath{\partial_{#2}{#1}}}
        }

\newcommand{\Fa}{For all }
\newcommand{\fa}{for all }

\newcommand{\fs}{for some }

\newcommand{\scth}{such that }

\newcommand\cf{\emph{cf.}~}

\newcommand\eg{\emph{e.g.}~}
\newcommand\ie{\emph{i.e.}~}

\newcommand\via{\emph{via}~}
\newcommand\loccit{\emph{loc.~cit.}}
\newcommand\opcit{\emph{op.~cit.}}

\newcommand\etc{\emph{etc.}}

\def\multi(#1,#2){\ifthenelse {\equal{#1}{0}}
                {{\mathbb Z}_2^{#2}}
                {\ifthenelse{\equal{#2}{0}}
                      {{\mathbb N}_0^{#1}}
                      {\ensuremath{{{\mathbb N}_0^{#1}\!\times\!{\mathbb Z}_2^{#2}}}}
                }
        }

\newcommand\vphi{\varphi}
\newcommand\vrho{\varrho}

\newcommand\vkappa{\varkappa}
\newcommand\eps{\varepsilon}
\newcommand\nats{\mathbb{N}}
\newcommand\ints{\mathbb{Z}}

\newcommand\reals{\mathbb{R}}
\newcommand\cplxs{\mathbb{C}}

\def\aff{{\mathbb A}}

\newcommand\sle{\leqslant}
\newcommand\sge{\geqslant}

\DMO\dom{\mathrm{dom}}
\DMO\rk{\mathrm{Rem}}
\DMO\Ad{\mathrm{Ad}}
\DMO\ad{\mathrm{ad}}
\DMO\GL{\mathrm{GL}}
\DMO\UGL{\mathrm{UGL}}
\DMO\USp{\mathrm{USp}}
\DMO\UU{\mathrm{U}}
\DMO\Sp{\mathrm{Sp}}
\DMO\SOSp{\mathrm{SOSp}}
\DMO\SO{\mathrm{SO}}
\DMO\id{\mathrm{id}}
\DMO\sdim{\mathrm{sdim}}
\DMO\sgn{\mathrm{sgn}}
\DMO\re{\mathrm{Re}}
\DMO\coker{\mathrm{coker}}
\DMO\im{\mathrm{im}}
\DMO\coim{\mathrm{coim}}
\DMO\codim{\mathrm{codim}}
\DMO\supp{\mathrm{supp}}
\DMO\str{\mathrm{str}}
\DMO\tr{\mathrm{tr}}
\DMO\vol{\mathrm{vol}}
\DMO\Spec{\mathrm{Spec}}
\DMO\res{\mathrm{res}}

\DMO\HC{\cat{HC}}
\DMO\CW{\cat{CW}}
\DMO\SMan{\cat{SMan}}
\DMO\Sets{\cat{Sets}}
\DMO\SSp{\cat{SSp}}
\DMO\SVec{\cat{SVec}}
\DMO\Shv{\cat{Shv}}

\DMO\WF{\mathrm{WF}}
\DMO\Op{\mathrm{Op}}

\newcommand{\lBr}{[\kern-.65ex[}
\newcommand{\rBr}{]\kern-.65ex]}

\makeatletter
\newcommand\Size[7][1]{
                                 \ifx#20%
                                        \def\r@l{}\def\r@m{}\def\r@r{}%
                                 \else%
                                    \ifx#21%
                                           \def\r@l{\bigl}\def\r@r{\bigr}\def\r@m{\bigm}%
                                    \else%
                                           \ifx#22%
                                                 \def\r@l{\Bigl}\def\r@r{\Bigr}\def\r@m{\Bigm}%
                                            \else%
                                                 \ifx#23%
                                                        \def\r@l{\biggl}\def\r@r{\biggr}\def\r@m{\biggm}%
                                                  \else
                                                        \ifx#24%
                                                        \def\r@l{\Biggl}\def\r@r{\Biggr}\def\r@m{\Biggm}%
                                                        \fi%
                                                  \fi%
                                            \fi%
                                      \fi%
                                 \fi%
                                 \ifx#10%
                                       \def\r@m{}%
                                 \fi%
                                 \r@l#3{#4}\r@m#5{#6}\r@r#7%
}%

\makeatother

\makeatletter
\def\Set@Scallop[#1]#2#3{{#1}\Parens{#2}{#3}}
\newcommand\DeclareScalableOperator[2]{%
  \expandafter\def\csname#1\endcsname{\@ifnextchar[{{#2}\Set@Scallop}{{#2}\Set@Scallop[{}]}}
}
\makeatother

\def\DSO{\DeclareScalableOperator}

\DSO{Presh}{\cat{Presh}}
\DSO{Sh}{\cat{Sh}}

\DSO{ShHom}{\sh Hom} 
\DSO{ShGHom}{\underline{\sh Hom}} 
\DSO{ShEnd}{\sh End} 
\DSO{ShGEnd}{\underline{\sh End}} 
\DSO{ShDer}{\sh Der} 
\DSO{ShGDer}{\underline{\sh Der}} 
\DSO{ShCliff}{\sh C\text{\emph{liff}}} 

\DSO{Ct}{\mathcal{C}} 
\DSO{Lp}{\mathbf{L}} 
\DSO{Hom}{\mathrm{Hom}} 
\DSO{End}{\mathrm{End}} 
\DSO{Aut}{\mathrm{Aut}} 
\DSO{Uenv}{\mathfrak{U}} 
\DSO{Cuenv}{\widehat{\mathfrak{U}}} 
\DSO{ABer}{\Abs0{\mathrm{Ber}}}
\DSO{Proj}{\mathrm{Proj}}
\DSO{Ber}{\mathrm{Ber}}
\DSO{Vol}{\mathrm{Vol}}
\DSO{Form}{\Omega}
\DSO{Ind}{\mathrm{Ind}}
\DSO{Coind}{\mathrm{Coind}}
\DSO{Der}{\mathrm{Der}}
\DSO{GDer}{\underline{\mathrm{Der}}}
\DSO{Alg}{\mathrm{Alg}}
\DSO{GHom}{\underline{\mathrm{Hom}}} 
\DSO{GEnd}{\underline{\mathrm{End}}} 
\DSO{Maps}{\mathrm{Maps}}
\DSO{Cliff}{\mathrm{Cliff}} 
\DSO{Jac}{\mathrm{Jac}} 
\DSO{Dist}{\mathcal{D}'}
\DSO{Sw}{\mathscr{S}}
\DSO{TDi}{\mathscr{S}'}

\newcommand\Set[3]{
                                 \Size{#1}{\{}{#2}{|}{#3}{\}}%
}%
\newcommand\Dual[3]{
                                 \Size[0]{#1}{\langle}{#2}{,}{#3}{\rangle}%
}%
\newcommand\Parens[2]{
  \Size[0]{#1}{(}{#2}{}{}{)}
}

\newcommand\Norm[2]{
  \Size[0]{#1}{\lVert}{#2}{}{}{\rVert}
}
\newcommand\Abs[2]{
  \Size[0]{#1}{\lvert}{#2}{}{}{\rvert}
}

\makeatletter

\makeatletter

\newif\if@smallmat
\newif\if@none
\newif\if@paren
\newif\if@brack
\newif\if@brace
\newif\if@vline
\newenvironment{Matrix}[2][1]
                                 {\ifx#20%
                                        \@smallmattrue%
                                  \else%
                                         \@smallmatfalse
                                  \fi%
                                  \ifx#11%
                                         \@nonefalse\@parentrue\@brackfalse\@bracefalse\@vlinefalse%
                                  \else%
                                       \ifx#12%
                                            \@nonefalse\@parenfalse\@bracktrue\@bracefalse\@vlinefalse%
                                        \else%
                                            \ifx#13%
                                                 \@nonefalse\@parenfalse\@brackfalse\@bracetrue\@vlinefalse%
                                            \else%
                                                 \ifx#14%
                                                       \@nonefalse\@parenfalse\@brackfalse\@bracefalse\@vlinetrue
                                                 \else%
                                                       \ifx#15%
                                                             \@nonefalse\@parenfalse\@brackfalse\@bracefalse\@vlinefalse%
                                                       \else%
                                                             \@nonetrue\@parenfalse\@brackfalse\@bracefalse\@vlinefalse%
                                                       \fi%
                                                 \fi%
                                            \fi%
                                        \fi%
                                   \fi%
                                   \if@smallmat%
                                        \if@none%
                                             \begin{smallmatrix}%
                                        \else%
                                            \if@paren%
                                                  \left(\begin{smallmatrix}%
                                            \else%
                                                  \if@brack%
                                                          \left[\begin{smallmatrix}%
                                                  \else%
                                                          \if@brace%
                                                               \left\{\begin{smallmatrix}%
                                                          \else%
                                                               \if@vline%
                                                                    \left\lvert\begin{smallmatrix}%
                                                                \else%
                                                                    \left\lVert\begin{smallmatrix}%
                                                                \fi%
                                                          \fi%
                                                  \fi%
                                            \fi%
                                        \fi%
                                   \else%
                                        \if@none%
                                             \begin{matrix}%
                                        \else%
                                            \if@paren%
                                                  \begin{pmatrix}%
                                            \else%
                                                  \if@brack%
                                                          \begin{bmatrix}%
                                                  \else%
                                                          \if@brace%
                                                               \begin{Bmatrix}%
                                                          \else%
                                                               \if@vline%
                                                                    \begin{vmatrix}%
                                                                \else%
                                                                    \begin{Vmatrix}%
                                                                \fi%
                                                          \fi%
                                                  \fi%
                                            \fi%
                                        \fi%
                                   \fi}%
                                  {\if@smallmat%
                                        \if@none%
                                             \end{smallmatrix}%
                                        \else%
                                            \if@paren%
                                                  \end{smallmatrix}\right)%
                                            \else%
                                                  \if@brack%
                                                          \end{smallmatrix}\right]%
                                                  \else%
                                                          \if@brace%
                                                               \end{smallmatrix}\right\}%
                                                          \else%
                                                               \if@vline%
                                                                    \end{smallmatrix}\right\rvert%
                                                                \else%
                                                                    \end{smallmatrix}\right\rVert%
                                                                \fi%
                                                          \fi%
                                                  \fi%
                                            \fi%
                                         \fi%
                                   \else%
                                        \if@none%
                                             \end{matrix}%
                                        \else%
                                            \if@paren%
                                                  \end{pmatrix}%
                                            \else%
                                                  \if@brack%
                                                          \end{bmatrix}%
                                                  \else%
                                                          \if@brace%
                                                               \end{Bmatrix}%
                                                          \else%
                                                               \if@vline%
                                                                    \end{vmatrix}%
                                                                \else%
                                                                    \end{Vmatrix}%
                                                                \fi%
                                                          \fi%
                                                  \fi%
                                            \fi%
                                        \fi%
                                   \fi}%
                         
\makeatother

\def\clap#1{\hbox to 0pt{\hss#1\hss}} 

\begin{document}

\title[Spherical superfunctions in rank one]{Asymptotics of spherical superfunctions on rank one Riemannian symmetric superspaces}

\author[Alldridge]
{Alexander Alldridge}

\address{Universit\"at zu K\"oln\\
Mathematisches Institut\\
Weyertal 86-90\\
50931 K\"oln\\
Germany}
\email{alldridg@math.uni-koeln.de}

\author[Palzer]
{Wolfgang Palzer}

\address{Universit\"at zu K\"oln\\
Mathematisches Institut\\
Weyertal 86-90\\
50931 K\"oln\\
Germany}
\email{palzer@math.uni-koeln.de}

\thanks{This research was supported by Deutsche Forschungsgemeinschaft, grant no.~DFG ZI 513/2-1 (AA), SFB/TR 12 ``Symmetries and University in Mesoscopic Systems'' (AA, WP), and by the institutional strategy of the University of Cologne within the German Excellence Initiative (AA)}

\subjclass[2010]{Primary 22E45, 17B15; Secondary 58A50}

\keywords{}

\begin{abstract}
  We compute the Harish-Chandra $c$-function for a generic class of rank-one purely non-compact Riemannian symmetric superspaces $X=G/K$ in terms of Euler $\Gamma$ functions, proving that it is meromorphic. Compared to the even case, the poles of the $c$-function are shifted into the right half-space. We derive the full asymptotic Harish-Chandra series expansion of the spherical superfunctions on $X$. In the case where the multiplicity of the simple root is an even negative number, they have a closed expression as Jacobi polynomials for an unusual choice of parameters. 
\end{abstract}

\maketitle

\section{Introduction}

Let $X=G/K$ be a Riemannian symmetric space. The mainstay of harmonic analysis on $X$ is the study of the asymptotic behaviour of the spherical functions $\phi_\lambda$ \cites{Hel84,Hel94,GV88}. These are $K$-invariant joint eigenfunctions of the $G$-invariant differential operators $D$ on $X$, \emph{viz.}
\[
  D\phi_\lambda=\Gamma(D)(\lambda)\phi_\lambda,
\]
where $\Gamma(D)(\lambda)$ is the value of the Harish-Chandra homomorphism. 

A basic observation is that $\phi_\lambda(a)$ admits an asympotic series expansion whose leading contribution is $c(\lambda)e^{(\lambda-\vrho)(a)}$, with $\vrho$ the half sum of the positive roots. Here, $c(\lambda)$ is Harish-Chandra's famous $c$-function. 

Moreover, $c(\lambda)$ admits an extension, as a meromorphic function, to the entire dual $\ger a^*$ of the complex Cartan subspace, due to the Gindikin--Karpelevi\v{c} formula
\[
  c(\lambda)=c_0\prod_\alpha\frac{2^{\lambda_\alpha-1}\Gamma(\lambda_\alpha)}{\Gamma\Parens1{\tfrac12(\lambda_\alpha+\frac{m_\alpha}2+1)}\Gamma\Parens1{\tfrac12(\lambda_\alpha+\frac{m_\alpha}2+m_{2\alpha})}},
\]
where the product is over all positive indivisible roots and $\lambda_\alpha\defi\frac{\Dual0\lambda\alpha}{\Dual0\alpha\alpha}$. This equation is basic in the proof of the fundamental theorems of harmonic analysis on Riemannian symmetric spaces: the support theorem for wave packets, the inversion, Paley--Wiener, Plancherel and Schwartz isomorphism theorems for the spherical Fourier transform. 

Although the situation is more complicated for the Helgason--Fourier transform (for $\tau$-spherical functions) and even more so for reductive non-Riemannian symmetric spaces, the basic philosophy of studying the asymptotics of spherical functions (or their replacements by more general Eisenstein integrals) remains valid, \cf Ref.~\cite{schlichtkrull-overview}. 

In the present paper, we study the asymptotics of spherical functions defined on a supersymmetric generalisation of Riemannian symmetric spaces. These play a basic role in the so-called Efetov supersymmetry method of condensed matter physics. Indeed, as shown by Zirnbauer \cites{Zir91,Zir92}, harmonic analysis on such superspaces can be used to give precise analytic expressions for the mean conductance of quasi-one dimensional disordered fermionic systems. 

Remarkably, contrary to the classical case, there are discrete contributions in the Plancherel decomposition, and this is visible and relevant in the physics of the corresponding systems. Since the Plancherel measure is governed by the $c$-function, the failure of absolute continuity can already be observed from the location of the $c$-function zeroes.

Indeed, let $(G,K)$ be one of the symmetric pairs listed in Table \ref{tab:rkone} below. Fix an indivisible restricted root $\alpha$ and $h_0$ with $\alpha(h_0)=1$. Let $\vrho=\frac12(m_\alpha\alpha+2m_{2\alpha}\alpha)$ where $m_\alpha$ and $m_{2\alpha}$ are the multiplicities of $\alpha$ and $2\alpha$, respectively, and identify $\lambda\equiv\lambda(h_0)$. Our first main result is the following theorem. 

\begin{ThA}
  For every $\Re\lambda>0$, the limit
  \[
     c(\lambda)=\lim_{t\longrightarrow\infty}\phi_\lambda(e^{th_0})e^{-t(\lambda-\vrho)}
  \] 
  exists. For some choice of $c_0\equiv c_0(\vrho)\neq0$, it is given by
  \begin{equation}\label{eq:clambda-intro}
    c(\lambda )=
    c_0\frac{2^{-\lambda }\Gamma (\lambda )}{\Gamma \left(\frac12\Parens1{\lambda +\frac{m_\alpha}2+1}\right)\Gamma \left(\frac12\Parens1{\lambda +\frac{m_\alpha}2+m_{2\alpha }}\right)}
  \end{equation}
  if $\alpha$ is anisotropic, and if $\alpha$ is isotropic, then it is given by 
  \begin{equation}\label{eq:clambda-intro-iso}
    c(\lambda)=c_0\lambda.
  \end{equation}
\end{ThA}

Formally, this result takes the same form as in the classical situation for the case of an anisotropic root $\alpha$. However, for an isotropic root $\alpha$ (of multiplicity $m_\alpha=-2$), one would expect $c(\lambda)\simeq(\lambda-1)$ from Equation \eqref{eq:clambda-intro} and the duplication formula. This differs from the true result in Equation \eqref{eq:clambda-intro-iso} by the absence of a `$\vrho$-shift'. This situation is similar for the Harish-Chandra homomorphism, see Ref.~\cite{a-hc}.

Moreover, in the present supersymmetric setting, $m_\alpha$ and hence $\vrho$ may be arbitrarily large negative numbers, drastically changing the asymptotic behaviour of $\phi_\lambda$. Moreover, due to the shift in the denominator in Equation \eqref{eq:clambda-intro}, $c(\lambda)$ picks up zeroes in the right half plane. These lead to discrete contributions in the Plancherel formula, as we will show in a forthcoming paper.

The class of symmetric pairs considered in Theorem A is a choice of rank one pairs that is generic in the sense that every rank one subpair of a reductive symmetric pair (of even type) associated with a choice of indivisible root (even, odd, ot both) is generated by copies of the pairs listed in Table \ref{tab:rkone} below: An anisotropic root $\alpha$ \scth $2\alpha$ is a root corresponds to the $\ger{gl}$ case; an anisotropic root \scth $2\alpha$ is not a root corresponds to the $\ger{osp}$ case, even it is purely odd (here the parameter $p=0$); finally, an isotropic root $\alpha$ corresponds to the $\ger{gl}(1|1)$ case. 

Choosing real forms of these symmetric pairs in such a way that the underlying symmetric spaces become Riemannian introduces extra conditions if one insists on taking real forms also of the odd part of the Lie superalgebra. However, these conditions are artificial from a physical point of view and moreover unnecessary for the analysis to go through. Therefore, we adopt the setting of `\emph{cs} manifolds' invented by Joseph Bernstein, that is, of real manifolds with complex sheaves of superfunctions. In many respects, this theory is parallel to that of real supermanifolds; however, there are some caveats, and we carefully lay the foundations to help the reader navigate these impasses in Section \ref{sec:pre}.

As the statement of Theorem A suggests, the general form of $c$-function for a supersymmetric symmetric pair is not given by a simple-minded generalisation of the Gindikin--Karpelevi\v c formula, since the contributions from isotropic roots are of a different form. Compare Ref.~\cite{ASc13} for details. 

The proof of Theorem A is somewhat more difficult than in the classical case. The reason lies again in the changed growth behaviour of the exponential $e^{t(\lambda-\vrho)}$: Upon parametrising the geodesic sphere at infinity of $X=G/K$ in stereographical coordinates, the integrals in question exhibit a singular behaviour. One therefore has to make a careful choice of cutoffs in a stereographical atlas parametrised by the Weyl group and keep track of this in the limit of $t\longrightarrow\infty$. 

In order for this to work, one needs to establish the Weyl group symmetry of the symmetric superfunctions $\phi_\lambda$ (\thmref{Cor}{spherical_inv}). Although this appears to be quite innocent, it requires the extension of the usual integral formul\ae{} for the Iwasawa and Bruhat decompositions to the `parameter-dependent' setting of supermanifolds over a general base supermanifold. For this reason, we are obliged to develop some foundational material in Subsections \ref{subs:rel-int} and \ref{subs:intfmla}. 

Once one has an explicit formula for the $c$-function, Harish-Chandra's series expansion of the spherical function carries over. (This also relies on \thmref{Cor}{spherical_inv}.) One has the following statement. 

\begin{ThB}
  Let $(G,K)$ be \scth the root $\alpha$ is even and $\Re\lambda>0$, $\lambda\notin\tfrac12\ints$. Then we have 
  \[
    \phi_\lambda|_{A^+}=\sum_{w\in W_0}\Phi_{w\lambda},\quad \Phi_{\pm\lambda}(e^{th_0})=e^{t(\pm\lambda-\vrho)}\sum_{\ell=0}^\infty\gamma_\ell(\lambda)e^{-2\ell t},
  \]
  where $A^+$ is the positive Weyl chamber, $W_0$ is the Weyl group, $\gamma_0(\pm\lambda)=c(\pm\lambda)$, and $\gamma_\ell(\lambda)$ follow a two-term recursion.
\end{ThB}

In most cases, $W_0=\{\pm1\}$, but for a suitable choice of parameters in case $G$ is an orthosymplectic supergroup, it can be trivial. A similar situation occurs when $\alpha$ is an odd anisotropic root (\thmref{Prop}{phi_iso}). In this case, the complexification of $G$ is $\GL(1|1,\cplxs)$. 

Remarkably, when $\vrho$ is a negative integer (as can happen in the $\ger{osp}$ case), the above series expansion is \emph{finite}, and $\phi_\lambda$ admits a simple closed expression in terms of Jacobi polynomials (\thmref{Cor}{sphfn-jacobi}) $P^{(a,b)}_n.$ Here, the parameters are chosen \scth $n=-\vrho$, so that in this case, the spherical superfunctions are exponential polynomials of a fixed degree independent of $\lambda$.

\medskip\noindent
We end this introduction by a synopsis of the paper's contents. Section \ref{sec:pre} is devoted to a brief collection of basic facts on \emph{cs} manifolds. We highlight  some points that to our knowledge are missing in the literature: \thmref{Prop}{mor-aff-hol} is a generalisation of Leites's theorem for morphisms from \emph{cs} manifolds to complex manifolds. Subsection \ref{subs:groups} contains an account of the exponential morphism for \emph{cs} Lie supergroups. Subsection \ref{subs:csforms} explains \emph{cs} forms of complex supergroups. Subsection \ref{subs:rel-int} is a brief account of the basic theory of Berezin fibre integrals. In Subsection \ref{subs:polar} a localization formula for Berezin integrals over superspheres is derived. 

Section \ref{sec:spherical} introduces the main players of this article, the spherical superfunctions. It begins by setting the stage for symmetric superpairs in Subsection \ref{subs:symmsup}. In Subsection \ref{subs:intfmla}, we collect some integral formul\ae{} necessary for the manipulation of the Harish-Chandra integral for $\phi_\lambda$. Notably, we generalise the classical formula for the Haar measure on a Lie group in exponential coordinates in \thmref{Prop}{pullback_nilpot}, the nilpotent case being of particular importance to us. The definition of the spherical functions is stated in Subsection \ref{subs:sphdef}; here, we also show a basic symmetry property of the spherical functions (\thmref{Cor}{spherical_inv}) that is important in the proof both of Theorem A (\thmref{Th}{gk-fmla}) and of Theorem B (\thmref{Th}{sph_sum}). 

Section \ref{sec:c_func} contains the statement and proof of Theorem A (\thmref{Th}{gk-fmla}). It is also the place where the main analytic problems have to be tackled and where the theory departs most from the classical cases (see the comments above). The theorem is stated in Subsection \ref{subs:gk-statement} and proved in Subsection \ref{subs:gk}. The proof proceeds case-by-case: the unitary case is contained in Subsubsection \ref{ssubs:c-unitary}, the orthosymplectic one in Subsubsection \ref{ssubs:osp}, and the case of an anisotropic purely odd root (the `$\GL(1|1,\cplxs)$ case') in Subsection \ref{ssubs:gl11}. The latter also contains the analytic expression for the spherical superfunctions $\phi_\lambda$ in this case. 

The final Section \ref{sec:hcseries} is devoted to the asymptotic series expansion of $\phi_\lambda$, that is, the proof of Theorem B (\thmref{Th}{sph_sum}). At this point, the main analytic difficulties have been overcome, so we can follow the ususal procedure of making a pertubation ansatz for the solutions of the eigenvalue equation for the radial part of the Laplacian on $X$, deducing a two-term recursion for the coefficients, and proving convergence \via Gangolli estimates. The estimates miraculously also go through in cases of negative multiplicity; if the half sum $\vrho$ of positive roots is a negative integer, the coefficients can even be computed explicitly, and the series terminates. In this case, we derive an expression for $\phi_\lambda$ in terms of Jacobi polynomials (\thmref{Cor}{sphfn-jacobi}) whose degree is fixed independent of $\lambda$.

\medskip\noindent
\emph{Acknowledgements.} This article is based in large parts on the second named author's doctoral thesis under the first named author's guidance. We wish to thank Martin Zirnbauer for constructive comments on early versions of our results and the anonymous referee for the careful reading of our paper and the detailed suggestions that helped to improve it.  

\section{Basic facts on \emph{cs} manifolds}\label{sec:pre}

In this section, we collect some basic facts on \emph{cs} manifolds and \emph{cs} Lie supergroups. As remarked above, this setting is necessary to cover the all the symmetric pairs we are interested in a Riemannian incarnation. We shall be suitably brief, including proofs only for those facts which so far have remained undocumented in the literature. 

\subsection{Basic definitions}

We will work with \emph{cs} manifolds and complex supermanifolds. The latter are covered by a wide literature, \eg by Refs.~\cites{ccf,manin,vsv}. The former, introduced by Joseph Bernstein, are covered to some degree in Ref.~\cite{deligne-morgan}. In many respects, they are similar to real supermanifolds, so one may follow the standard texts on that subject \cites{deligne-morgan,leites,vsv}. The differences that exist are quite subtle, and we will comment on these to the extent required in our applications. Specifically, we shall follow the conventions of Ref.~\cite{ahw-sing}.

In particular, we will work in the framework of $\cplxs$-superspaces, that is, of locally superringed spaces over $\Spec\cplxs$. Our notation will always be $X=(X_0,\sh O_X)$ for $\cplxs$-superspaces, and $\vphi=(\vphi_0,\smash{\vphi^\sharp}):X\longrightarrow Y$ for morphisms, where we may consider $\vphi^\sharp$ as a sheaf map $\sh O_Y\longrightarrow \vphi_{0*}\sh O_X$ or $\vphi_0^{-1}\sh O_Y\to O_X$ by the fundamental adjunction of the direct and inverse image functors, see \eg Refs.~\cite{bredon}*{I.4, Equation (5)} or \cite{iversen}*{II.4, Theorem 4.8}.

Somewhat abusing notation, we will also write $X_0$ for the reduced superspace associated with $X$, see \cite{ahw-sing}*{Construction 3.9}. There is a canonical closed embedding 
\[
  j_{X_0}:X_0\longrightarrow X
\]
and for any $f\in\sh O_X(U)$, where $U\subseteq X_0$ is open, and any $x\in U$, we denote by $f(x)\defi \smash{j_{X_0}^\sharp(f)(x)}\in\vkappa(x)\defi\sh O_{X,x}/\ger m_{X,x}$ the \Define{value} of $f$ at $x$. In the cases of interest to us, we will always have $\vkappa(x)=\cplxs$.

We define the \emph{cs} and complex affine superspaces as 
\[
  \aff^{p|q}\defi(\reals^p,\sh C^\infty_{\reals^p}\otimes\textstyle\bigwedge(\cplxs^q)^*),\quad
  \aff_{hol}^{p|q}\defi(\cplxs^p,\sh H_{\cplxs^p}\otimes\textstyle\bigwedge(\cplxs^q)^*),
\] 
where $\sh C^\infty_{\reals^p}$ is the sheaf of smooth complex-valued functions on $\reals^p$ and $\sh H_{\cplxs^p}$ is the sheaf of holomorphic functions on $\cplxs^p$, and we consider the Euclidean topology. Then, by definition, a \Define{\emph{cs} manifold} is a $\cplxs$-superspace locally isomorphic to $\aff^{p|q},$ whereas a \Define{complex supermanifold} is one locally isomorphic to $\smash{\aff^{p|q}_{hol}}.$ 

We are tempted speak of \emph{cs} manifolds simply as ``supermanifolds'', but for the sake of convention, we will stick in this paper to the original appellation, which was introduced by Joseph Bernstein, \cf Ref.~\cite{deligne-morgan}.

More generally, given two \emph{cs} manifolds $X$ and $S$, a \Define{\emph{cs} manifold over $S$} is a morphism $X\to S$ that is locally in $X$ isomorphic to the projection of a direct product $S\times Y\to S$. In this case $\dim_SX\defi\dim Y$ is called the \Define{fibre dimension}. Usually, we will just write $X/S$ without mentioning the structural morphism explicitly, denoting it by $p_X$ where necessary. In particular, $\smash{\aff^{p|q}_S}\defi S\times\smash{\aff^{p|q}}$ is a \emph{cs} manifold over $S$.  A \Define{morphism $f:X/S\to Y/S$ over $S$} is a morphism $f:X\to Y$ respecting the structural morphisms. Given an open embedding $X/S\to\smash{\aff^{p|q}_S}/S$ over $S$, the pullback of the standard coordinate functions on $\aff^{p|q}$ is called a \Define{fibre coordinate system} (over $S$). We recover the usual \emph{cs} manifolds upon setting $S=*$, the terminal object of the category of $\cplxs$-superspaces (\ie the singleton space, together with the constant sheaf $\cplxs$), but the relative point of will be important in Subsections \ref{subs:tanfun}, \ref{subs:rel-int}, \ref{subs:groups}, and \ref{subs:intfmla}.

In Ref.~\cite{leites}, the corresponding notion for real supermanifolds is discussed under the name of \emph{families}. In Ref.~\cite{ahw-sing}, which we follow, it is defined over base superspaces $S$ more general than \emph{cs} manifolds. 

A coordinate-free view on affine superspaces will be useful. To that end, we define a \emph{cs} vector space to be a real super-vector space $V=V_\ev\oplus V_\odd$ with a distinguished complex structure on $V_\odd$. We set 
\[
  \aff(V)\defi\Parens1{V_\ev,\sh C^\infty_{V_\ev}\otimes\textstyle\bigwedge V_\odd^*}
\]
for any finite-dimensional \emph{cs} vector space $V$. In particular, $\aff^{p|q}=\aff(\reals^p\oplus\Pi\cplxs^q)$. Similarly, for any finite-dimensional complex super-vector space $V,$ we define
\[
  \aff_{hol}(V)\defi\Parens1{V_\ev,\sh H_{V_\ev}\otimes\textstyle\bigwedge V_\odd^*},
\]
so that $\aff^{p|q}_{hol}=\aff_{hol}(\cplxs^{p|q})$, where we write $\cplxs^{p|q}\defi\cplxs^p\oplus\Pi\cplxs^q$.

We take note of the following trivial observation
\[
  V_\ev=0\ \Longrightarrow\ \aff(V)=\aff_{hol}(V).
\]

\subsection{Morphisms to affine superspaces}\label{subs:morphisms}

One essential respect in which \emph{cs} manifolds differ from real or complex supermanifolds is that morphisms to affine superspace are characterised in a subtly different fashion. Namely, for any $\cplxs$-superspace $X$, denote by $\sh O_{X,\reals}$ the sheaf whose local sections are the superfunctions $f$ (\ie local sections of $\sh O_X$) whose values $f(x)$ are real, for any $x\in X_0$ at which $f$ is defined. Then we have the following generalisation of Leites's theorem on morphisms.

\begin{Prop}[mor-aff][\protect{\cite{ahw-sing}*{Lemma 3.13, Corollary 4.15}}]
  Let $T$ be a \emph{cs} manifold and denote the standard coordinate system on $\aff^{p|q}$ by $(t_a)$. Then the natural map
  \[
    \aff^{p|q}(T)\defi\Hom0{T,\aff^{p|q}}\longrightarrow\Gamma\Parens1{\sh O_{X,\reals,\ev}^p\times\sh O_{X,\odd}^q}:f\longmapsto(f^\sharp(t_a))
  \]
  is a bijection. Here and in what follows, $\mathrm{Hom}$ denotes morphisms of $\cplxs$-superspaces and $\Gamma$ denotes the global sections functor for sheaves. 
\end{Prop}

Expanding on the above conventions, we let, for any \emph{cs} manifolds $X$ and $T$, $X(T)\defi\Hom0{T,X}$ denote the set of morphisms $T\longrightarrow X$. We call any such morphism $x:T\longrightarrow X$ a $T$-valued point and write $x\in_TX$. For a morphism $f:X\to Y$ and $x\in_TX$, we denote by $f(x)\defi f\circ x\in_TY$. Given a $T$-valued point $x\in_TX$ where $T/S$ is a \emph{cs} manifold over $S$, we write $x\in_{T/S}X/S$ (or $x\in_{T/S}X$ for brevity) if $x$ is over $S$ as a morphism $T\to X$. We will use this suggestive notation constantly. Together with the Yoneda lemma, it provides a convenient way of defining and handling morphisms. We will be constantly using this point of view, so \thmref{Prop}{mor-aff} is basic to our study. 

A coordinate-free version of \thmref{Prop}{mor-aff} can be formulated as follows: Following Ref.~\cite{ahw-sing}, we define  
\[
  (V\otimes W)_\ev\defi V_\ev\otimes_\reals W_\ev\oplus V_\odd\otimes_\cplxs W_\odd.
\]
With this notation, the above proposition admits the following reformulation.

\begin{Cor}[mor-aff-coordfree][\protect{\cite{ahw-sing}*{Corollary 3.24}}]
  Let $V$ be a finite-dimensional \emph{cs} vector space and $T$ a \emph{cs} manifold. Then there is a natural bijection 
  \[
    \aff(V)(T)=\Hom0{T,\aff(V)}\longrightarrow\Gamma\Parens1{(\sh O_{T,\reals}\otimes V)_\ev}.
  \]
\end{Cor}

On the other hand, as the holomorphic version of Leites's theorem states, morphisms from a \emph{complex} supermanifold $T$ to $\smash{\aff^{p|q}_{hol}}$ constitute the same data as tuples of superfunctions on $T$ \emph{without any reality condition}. Remarkably, this carries over to the situation where $T$ is a \emph{cs} manifold, by the results of Ref.~\cite{ahw-sing}. 

\begin{Prop}[mor-aff-hol]
  Let $V$ be a finite-dimensional complex super-vector space and $T$ a \emph{cs} manifold. Then there is a natural bijection 
  \[
    \aff_{hol}(V)(T)=\Hom0{T,\aff_{hol}(V)}\longrightarrow\Gamma\Parens1{(\sh O_T\otimes_\cplxs V)_\ev}.
  \]
\end{Prop}

\begin{proof}
  Broadly following Ref.~\cite{ahw-sing}, a $\cplxs$-superspace $T$ is called \emph{holomorphically regular} if the above natural map is bijective for any $V$ and any open subspace of $T$. 

  Thus, our claim is that any \emph{cs} manifold is holomorphically regular. We will use the terminology of Ref.~\cite{ahw-sing}. Since its use will be localized to this proof, we refer to that article for any undefined notions. Firstly, we may assume that $V=\cplxs^n,$ by Lemma 4.15 (\opcit). Then $\smash{\aff_{hol}(V)=\aff^n_{hol}}$ and $\smash{\aff^n_{hol}(T)=\aff^n_{hol}(T^\ev)},$ where $T^\ev\defi(T_0,\sh O_{T,\ev})$ is the even part of $T$, so it is sufficient to prove that $T^\ev$ is holo\-morph\-ical\-ly regular. 

  But by Proposition 5.24 (\opcit), $T^\ev$ (and hence, any open subspace thereof) admits a tidy embedding into some $\aff^N.$ In view of Proposition 5.11 and Lemma 5.2 (\opcit), it is therefore sufficient to prove that $\aff^N$ is holomorphically regular. By Lemma 5.2 (\opcit) again, the natural map is injective for any open subspace of $\aff^N$. Let us prove that it is also surjective.

  So, let $U\subseteq\reals^N$ be open and $f_1,\dotsc,f_n\in\sh O_{\aff^N}(U)=\Ct[^\infty]0{U,\cplxs}.$ We may define $f_0\defi(f_1,\dotsc,f_n):U\to\cplxs^n,$ and this map is smooth. Setting
  \[
    f^\sharp_V(h)\defi h\circ\Parens1{f_0|_{\smash{f_0^{-1}(W)}}},\quad h\in\sh O_{\aff^n_{hol}}(W)=\sh H(W),
  \]
  for any open $W\subseteq\cplxs^n,$ defines a morphism $f=(f_0,f^\sharp):\aff^N|_U\to\aff^n_{hol}$ \scth $f^\sharp(z_a)=f_a,$ where $(z_a)$ are the standard coordinates on $\aff^n_{hol}.$
\end{proof}

For any finite-dimensional complex super-vector space $V$, we define a set-valued cofunctor $\aff^\cplxs(V)$ on the category of \emph{cs} manifolds by setting 
\[
  \aff^\cplxs(V)(T)\defi\Gamma\Parens1{(\sh O_T\otimes_\cplxs V)_\ev}
\]
on objects $T$, and by the obvious definition on morphisms. In these terms, \thmref{Prop}{mor-aff-hol} states that $\aff^\cplxs(V)$ is the restriction to the category of \emph{cs} manifolds of the point functor of the complex supermanifold $\aff_{hol}(V)$. (Of course, the restriction of the point functor of the \emph{algebraic} affine superspace $\Spec S(V^*)$ to \emph{cs} manifolds is also the same, but this is somewhat less remarkable.)

\thmref{Prop}{mor-aff-hol} has the following consequence, which will be important in the applications to supergroups below. 

\begin{Cor}[complex-prod]
  Let $X$ and $Y$ be complex supermanifolds and $T$ a \emph{cs} manifold. Then the natural map
  \[
    (X\times Y)(T)\longrightarrow X(T)\times Y(T)
  \]
  is a bijection. 
\end{Cor}

\subsection{The tangent functor}\label{subs:tanfun}

Preparing for our discusson of the exponential morphism of \emph{cs} Lie supergroups below, we introduce the `total space' of the `tangent bundle'. A subtlety for \emph{cs} manifolds $X$ is that such a `total space' representing the (complex) tangent sheaf does not literally exist as a \emph{cs} manifold. In this section, we introduce a replacement on the level of functors. We then single out a representable subfunctor, the `real tangent bundle'. We will work with \emph{cs} manifolds over an arbitrary base $S$. This will be used to the define the exponential morphism. It will also be important for the machinery of Berezinian fibre integration that we introduce in Subsection \ref{subs:rel-int}, on which our proof of the Weyl symmetry of the spherical superfunctions in \thmref{Cor}{spherical_inv} hinges. 

\begin{Def}[tanfun][the tangent functor]
  Let $X/S$ be a \emph{cs} manifold over $S$. We let $\sh T_{X/S}$ be the sheaf of $p_{X,0}^{-1}\sh O_S$-linear superderivations of $\sh O_X$, and call this the \Define{tangent sheaf over $S$}. Locally, it is spanned by fibre coordinate derivations and hence locally free. We define the \Define{tangent functor} $T_SX$ of $X$ over $S$ as the set-valued cofunctor on \emph{cs} manifolds given on objects $T/S$ by 
  \[
    (T_SX)(T/S)\defi\Set1{(x,v)}{x\in_{T/S}X,v\in\Gamma\Parens1{(x^*\sh T_{X/S})_\ev}}
  \]
  and by the obvious definition on morphisms. There are canonical morphisms of functors
  \[
    0_X:X\longrightarrow T_SX:x\longmapsto(x,0),\quad\pi_{X/S}:T_SX\longrightarrow X:(x,v)\longmapsto x,
  \]
  called the \Define{zero section} and the \Define{projection}, respectively. If $S=*$, we omit the corresponding subscripts. 

  The construction of $T_SX$ is functorial: For any morphism $\vphi:X/S\to Y/S$ of \emph{cs} manifolds over $S$, the tangent morphism $T_S\vphi:T_SX\to T_SY$ is defined by 
  \[
    (T_S\vphi)_T(x,v)\defi\Parens1{\vphi(x),T_{x/S}\vphi(v)},\quad (T_{x/S}\vphi)(v)(f)\defi v\Parens1{\vphi^\sharp(f)}
  \]
  \fa $f\in(\vphi_0^{-1}\sh O_Y)(U)$, $U\subseteq Y_0$ open. Here, we recall that $x^*\sh T_{X/S}$ is the sheaf of vector fields over $S$ and along $x$.

  Clearly, the functor $T_S:X\mapsto T_SX$ preserves fibre products over $S$, so $T_SS=S$ and $T*=*$, the point functor of the singleton space.
\end{Def}

This definition of the tangent functor is compatible with the definition of the tangent spaces. If $x\in X_0$ is a point, considered as a morphism $*\longrightarrow X$, then the fibre product of functors $(T_SX)_x=*\times_XT_SX$ is given by 
\begin{equation}\label{eq:tangent-comp}
  \begin{split}
    (T_SX)_x(T)&=(*\times_XT_SX)(T)=\Gamma\Parens1{(x_T^*\sh T_{X/S})_\ev}\\
    &=\Set1{\vphi\in\Hom[_S]0{T\times\Spec\mathbb D_\ev,X}}{\vphi|_{\eps=0}=x_T}\\
    &=\aff^\cplxs(T_{x/S}X)(T).  
  \end{split}    
\end{equation}
Here, $x_T$ is the specialization of $x$, \ie the composite $T\longrightarrow *\longrightarrow X$, and $T_{x/S}X$ is the tangent space of $X$ over $S$ at $x$---that is, the super-vector space over $\cplxs$ whose homogeneous elements are the $\sh O_{S,p_0(x)}$-linear maps $v:\sh O_{X,x}\to\cplxs$ \scth 
\[
  v(fg)=v(f)g+(-1)^{\Abs0f\Abs0v}fv(g),\quad\forall f,g.
\]
In other words, the fibre of $T_SX\longrightarrow X$ over $x$ is the $\cplxs$-affine superspace of $T_{x/S}X$. Moreover, $\Spec\mathbb D_\ev=(*,\cplxs[\eps]/(\eps^2))$, and $(\cdot)|_{\eps=0}$ denotes the canonical morphism $T\longrightarrow T\times\Spec\mathbb D_\ev$ that pulls back functions in $\sh O_T[\eps]/(\eps^2)$ by dropping the linear term in $\eps$. The last of the equalities in Equation \eqref{eq:tangent-comp} amounts to the routine check that the map $\smash{\vphi^\sharp}=x_T^\sharp+\eps v$ is an algebra morphism if and only if $v$ is a vector field over $S$ along $x_T$.

Similarly, one checks that for any morphism $\vphi:X\to Y$, the action of $T\vphi$ on the fibre $(T_SX)_x$ is simply given by the application of the functor $\aff^\cplxs(\cdot)$ to the ordinary tangent map $T_{x/S}\vphi:T_xX\to T_{\vphi_0(y)}Y$.

\medskip
Besides the (complex) tangent space $T_{x/S}X$ over $S$, we consider the \Define{real tangent space}, defined as the \emph{cs} vector space
\[
  T_{x/S}^\reals X\defi\Set1{v\in T_{x/S}X}{v(\sh O_{X,\reals,\ev})\subseteq\reals}
\]
of real tangent vectors. More generally, we introduce the following.

\begin{Def}[real-tan][the real tangent bundle]
  The \Define{real tangent sheaf} $\smash{\sh T_{X/S}^\reals}$ of $X$ over $S$ as the subsheaf of $\sh T_{X/S}$ whose local sections over an open set $U\subseteq X_0$ are those $v\in\sh T_{X/S}(U)$ \scth 
  \[
    v(\sh O_{X,\reals,\ev}(V))\subseteq\sh O_{X,\reals}(V)
  \]
  \fa open subsets $V\subseteq U$. 

  If $(x_a)=(u,\xi)$ are local fibre coordinates of $X$ over $S$ defined on $U$, then 
  \[
    \sh T^\reals_{X/S}\big|_U=\bigoplus_i\sh O_{X,\reals}\frac\partial{\partial u_i}\oplus\bigoplus_j\sh O_X\frac\partial{\partial\xi_j},
  \]
  where $\smash{\tfrac\partial{\partial u_i}}$ and $\smash{\tfrac\partial{\partial\xi_j}}$ are the fibre coordinate derivations introduced by $(u,\xi)$. In this sense, $\smash{\sh T^\reals_{X/S}}$ is a `locally free graded module over the ring extension $(\sh O_X,\sh O_{X,\reals})$'.

  We define the \Define{real tangent functor} $T^\reals_SX$ as the subfunctor of $T_SX$ given on objects by 
  \[
      (T^\reals_SX)(T/S)\defi\Set1{(x,v)\in T_SX(T/S)}{v\in\Gamma\Parens1{(x^*\sh T^\reals_{X/S})_\ev}}.
  \]
  The condition on $v$ in the last equation amounts to 
  \[
    v\Parens1{(x_0^{-1}\sh O_{X,\reals,\ev})(U)}\subseteq\sh O_{T,\reals}(U)
  \]
  \fa open subsets $U\subseteq T_0$. Due to the local freeness of $\smash{\sh T^\reals_{X/S}}$, the functor $T^\reals_SX$ is representable by a \emph{cs} manifold, which is the total space of a vector bundle over $X$. This vector bundle, also denoted by $T^\reals_SX$, is called the \Define{real tangent bundle of $X.$} 

  By \thmref{Prop}{mor-aff} or \thmref{Cor}{mor-aff-coordfree}, its fibre at $x\in X_0$ is computed to be 
  \[
    (T_S^\reals X)_x=\aff\Parens1{T^\reals_{x/S}X},
  \]
  the \emph{cs} affine superspace associated with the \emph{cs} vector space $T_{x/S}^\reals X$.  
\end{Def}

\subsection{The exponential morphism of a \emph{cs} Lie supergroup}\label{subs:groups}

It is a well-known fact that there is an equivalence between the categories of Lie supergroups and supergroup pairs, \cf Ref.~\cites{ccf,deligne-morgan,kostant,koszul}. This is true in the real and complex, smooth and analytical cases, and the arguments valid in these cases carry over unchanged to the setting of \emph{cs} Lie supergroups. 

However, for applications, in particular those to Harmonic Superanalysis we are concerned with in this paper, it useful to have a `geometric' view on this equivalence. This uses the exponential morphism and the Campbell--Hausdorff series, neither of which has as yet been given a treatment in the literature in this setting. 

The exponential morphism will also be used extensively in Section \ref{sec:c_func}, for the derivation of the Iwawasa $H$-projection and the integration in exponential coordinates on the nilpotent Iwawasa $N$ supergroup (on the basis of \thmref{Prop}{pullback_nilpot}), both of which are primordial for the $c$-function asymptotics. Moreover, the \emph{cs} Lie supergroups we will consider (see below for the definition) will be `forms' of complex Lie supergroups. We will widely use their functors of points, and the avalability of a manageable expression thereof will be crucial. As we show below, in \thmref{Prop}{csform-points}, such an expression can be derived by means of the exponential morphism. We begin with some basic definitions. 

\begin{Def}
  By definition, a \Define{\emph{cs} Lie supergroup} is a group object $G$ in the category of \emph{cs} manifolds. 
\end{Def}

Let $G$ be a \emph{cs} Lie supergroup. We set 
\[
  \ger g\defi T_1G,\quad\ger g_\reals\defi T_1^\reals G,
\]
where we write $T^\reals_1G=T^\reals_{1/*}G$ and otherwise use the notation from Subsection \ref{subs:tanfun}. We have 
\[
  \ger g=\Set1{x}{g=1+\eps x_\ev+\tau x_\odd\in G(\Spec\mathbb D)},\quad\mathbb D\defi\cplxs[\eps|\tau]/(\eps^2,\eps\tau),
\]
where we set $\Spec\mathbb D\defi (*,\mathbb D)$. Hence, for $x,y\in\ger g$, we may define $[x,y]$ by 
\[
    a=1+\eps[x,y]_\ev+\tau[x,y]_\odd,\quad a\defi(ghg^{-1}h^{-1})|_{\eps=\eps_1\eps_2=\tau_1\tau_2,\tau=\tau_1\eps_2=\eps_1\tau_2},
\]
where
\[
  g\defi1+\eps_1x_\ev+\tau_1x_\odd,\quad h\defi1+\eps_2y_\ev+\tau_2y_\odd.
\]
With this bracket, $\ger g$ acquires the structure of a complex Lie superalgebra. Moreover, we have $\ger g_{\reals,\ev}=\ger g_0$ and $\ger g_\ev=\ger g_0\otimes_\reals\cplxs$, where $\ger g_0$ is the Lie algebra of the real Lie group $G_0$. Note that by \thmref{Cor}{mor-aff-coordfree}, we have 
\[
  \aff(\ger g_\reals)(\aff^{0|q})=T_1^\reals((\Pi T^q)G)_0,
\]
the Lie algebra of the real Lie group $((\Pi T)^qG)_0=G(\aff^{0|q})$.

When $\vphi:G\to H$ is a morphism of \emph{cs} Lie supergroups, then $\vphi(1)=1$, so that we have a map $d\vphi\defi T_1\vphi:\ger g=T_1G\to T_1H=\ger h$. Using the definition of the bracket, it is immediate that $d\vphi$ is a morphism of complex Lie superalgebras. Moreover, we have $d\vphi(\ger g_\reals)\subseteq\ger h_\reals$.

We may now transfer the definition of the exponential morphism of a real Lie supergroup \cite{GW12} to the setting of \emph{cs} Lie supergroups. We shall be terse on the parts that are similar, highlighting only the use of the tangent functor. 

\begin{Prop}[exp]
  Let $G$ be a \emph{cs} Lie supergroup. There is a unique morphism 
  \[
    \exp_G:\aff(\ger g_\reals)\to G
  \]
  of \emph{cs} manifolds whose action on $\aff^{0|q}$-valued points is precisely the exponential map of the real Lie group $((\Pi T)^qG)_0=G(\aff^{0|q})$. The morphism $\exp_G$ is a local isomorphism in a neighbourhood of zero.
\end{Prop}

\begin{proof}
  The uniqueness is obvious, since the $\aff^{0|q}$ form a set of generators for the category of \emph{cs} manifolds. 

  For the existence, we construct a certain even vector field on $G\times\aff(\ger g_\reals)$. If $m$ is the multiplication of $G$, we consider
  \[
    0_G\times\iota:G\times\aff^\cplxs(\ger g)\to TG\times TG,\quad Tm:TG\times TG=T(G\times G)\to TG,
  \]
  where $\iota$ is the canonical morphism $\aff^\cplxs(\ger g)=(TG)_1\to TG$. 

  Let $L$ denote their composite. Explicitly, on points, it is given by 
  \[
    L(g,v)=Tm(g,1,0,v)=\Parens1{g,v_g},\quad v_g(f)\defi({\id}\otimes v)\Parens1{m(g,\cdot)^\sharp(f)}.
  \]
  The inverse morphism $L^{-1}:TG\to G\times\aff^\cplxs(\ger g)$ is given on points by 
  \[
    L^{-1}(g,v)=(g,v_{g^{-1}}),
  \]
  so that $L$ is an isomorphism of functors. 

  Now, we define the morphism $\Theta:G\times\aff(\ger g_\reals)\longrightarrow TG$ as the composite 
  \begin{center}
    \begin{tikzcd}
      G\times\aff(\ger g_\reals)\rar{}&G\times\aff(\ger g)\rar{L}&TG.
    \end{tikzcd}
  \end{center}
  By the Yoneda Lemma, $\Theta$ corresponds to a unique element $(x,v)$ of 
  \[
    T(G)(G\times\aff(\ger g_\reals)).
  \]
  It is clear that $\pi_G\circ\Theta=p_1:G\times\aff(\ger g_\reals)\to G$, so that $v\in\Gamma\Parens1{(p_1^*\sh T_G)_\ev}$. 

  We may promote $v$ to an even vector field on $G\times\aff(\ger g_\reals)$ over $\aff(\ger g_\reals)$. Arguing as in \cite{GW12}*{Lemma 4.1}, one sees that $v$ is real and complete. Let 
  \[
    \gamma_v:\aff^1\times G\times\aff(\ger g_\reals)\to G\times\aff(\ger g_\reals)
  \]
  be its global flow. We define $\exp_G$ as the composite
  \[
    \begin{tikzcd}[column sep=large]
      \aff(\ger g_\reals)\rar{(1,1)\times{\id}}&\aff^1\times G\times\aff(\ger g_\reals)\rar{\gamma_v}&G\times\aff(\ger g_\reals)\rar{p_1}& G.
    \end{tikzcd}    
  \]
  Using the naturality of the morphism $L$, we see that $v$ corresponds on $G(\aff^{0|q})$ to a similarly defined vector field that is known to generate the exponential flow \cite{kms}. This implies the claim. 
\end{proof}

\begin{Cor}[exp-comm]
  Let $\vphi:G\to H$ be a morphism of \emph{cs} Lie supergroups. Then 
  \[
    {\exp_H}\circ\aff(d\vphi)=\vphi\circ\exp_G.
  \]
\end{Cor}

\begin{proof}
  It is sufficient to check the equality on $\aff^{0|n}$-valued points, where it is immediate from \thmref{Prop}{exp} and the classical theory \cite{hn}. 
\end{proof}

The following is a form of the equivalence of supergroups and supergroup pairs, stating that \emph{any supergroup $G$ is $G_0$-equivariantly split}. (In fact, it gives the essential bijection on objects; the equivalence on morphisms then follows from \thmref{Cor}{exp-comm}.)

\begin{Cor}[grp_pnt]
  Let $G$ be a \emph{cs} Lie supergroup. Then the morphism
  \[
    G_0\times\aff(\ger g_\odd)\longrightarrow G:(g,x)\longmapsto g\exp_G x
  \]
  is a $G_0$-equivariant isomorphism of \emph{cs} manifolds. 
\end{Cor}

\begin{proof}
  Since $d\exp_G=\id_\ger g$ by \thmref{Prop}{exp} and the morphism is $G_0$-equivariant, it has invertible tangent map at every point. The underlying map is ${\id}_{G_0}$, so the claim follows from the inverse function theorem \cites{ccf,deligne-morgan,leites}.
\end{proof}

We end this subsection by a discussion of the \Define{adjoint action}. This is the action $\Ad_G\defi\Ad$ of $G$ on $\aff^\cplxs(\ger g)=(TG)_1$, given as follows: For $g\in_SG$ and $x\in(TG)_1(S)$, regard $x$ as an element of $G(S\times\Spec\mathbb D_\ev)$ \scth $x|_{\eps=0}=1_S$. (See the comments after Equation \eqref{eq:tangent-comp} for explanation.) Then 
\[
  \Ad(g)(x)\defi gxg^{-1}\in G(S\times\Spec\mathbb D_\ev),
\]
where $g$ is specialized to $S\times\Spec\mathbb D_\ev$ \via the first projection. It follows that 
\[
  \Ad(g)(x)|_{\eps=0}=gg^{-1}=1_S, 
\]
so $\Ad(g)(x)\in_S(TG)_1=\aff^\cplxs(\ger g)$. Using the definitions, it is easy to check that for any morphism $\vphi:G\to H$ of \emph{cs} Lie supergroups, we have 
\[
  \Ad_H(\vphi(g))(\aff^\cplxs(d\vphi)(x)))=\aff^\cplxs(d\vphi)\Parens1{\Ad_G(g)(x)}
\]
for any $g\in_SG$, $x\in_S\aff^\cplxs(\ger g)$.

If we let $G_0\times G_0$ act on $G_0\times\aff(\ger g_\odd)$ by $(g_1,g_2)(g,x)=(g_1gg_2^{-1},\Ad(g_2)(x))$ for $(g_1,g_2)\in_TG_0\times G_0$, then we see that the isomorphism stated in \thmref{Cor}{grp_pnt} is even $(G_0\times G_0)$-equivariant. 

Remarkably, the adjoint action of $G$ passes to $\aff(\ger g_\reals)$. Namely, let $g\in_SG$ and $x\in_S\aff(\ger g_\reals)$. We need to check that for any $f\in\Gamma(\sh O_{G,\reals,\ev})$, we have $y(f)\in\Gamma(\sh O_{S,\reals})$ for $y\defi\Ad(g)(x)$. By the definition, we have
\[
  y(f)=({\id}\otimes x)(f(g(\cdot)g^{-1})),
\]
where ${\id}\otimes x$ denotes the promotion of $x$ to a vector field over $S$ (which exists because of the local freeness of the tangent sheaf over $S$, \cf \cite{leites}*{Lemma 2.2.3}).  

Since $f(g(\cdot)g^{-1})=(m(g,m(\cdot,g^{-1})))^\sharp(f)$ is an even and real-valued superfunction on $S\times G$, the assertion is immediate from the assumption on $x$.

\begin{Cor}[exp-ad]
  Let $G$ be a \emph{cs} Lie supergroup. If $g\in_SG$ and $x\in_S\aff(\ger g_\reals)$, then 
  \[
    \exp_G\Ad(g)(x)=g(\exp_Gx)g^{-1}.
  \]
\end{Cor}

\begin{proof}
  By the above, both sides of the equation are well-defined. Hence, it is sufficient to check the equality on $\aff^{0|q}$-valued points, for every $q$. But on that level, it follows from \thmref{Prop}{exp} and the classical theory \cite{hn}. Alternatively, one may use a parameter version of \thmref{Cor}{exp-comm}.
\end{proof}

\begin{Rem}
  Although $\ger g_\reals$ is in general not a Lie algebra, the bracket is well-defined on $\aff(\ger g_\reals)$. Indeed, $\ger g$ has a homogeneous basis $(e^a)$ contained in $\ger g_\reals$. Then $x,y\in\Gamma\Parens1{\sh O_{T,\odd}\otimes\ger g_\odd}$ admit representations $x=\sum_ax_ae^a$ and $y=\sum_ay_ae^a$, and 
  \[
    [x,y]=-\sum_{abc}x_ay_bC^{ab}_ce^c  
  \]
  where $C^{ab}_c\in\cplxs$ are the structure constants of $\ger g$. But $x_ay_bC^{ab}_c$ has the value zero, so that it is a section of $\sh O_{T,\reals,\ev}$.
\end{Rem}

\subsection{Complex supergroups and \emph{cs} forms}\label{subs:csforms}

Our interest in \emph{cs} supergroups comes from the fact that their are many `\emph{cs} forms' of complex Lie supergroups, whereas there are comparatively few real forms. (See the examples we will be considering in Section \ref{sec:c_func} and beyond, and compare Ref.~\cite{deligne-morgan} for a more extensive list of reasons. In any case, any real Lie supergroup defines a \emph{cs} Lie supergroup by the complexification of its structure sheaf.)

Since we will encounter many such \emph{cs} forms, it will be useful to have a uniform description of their functors of points. We will derive such a description by the use of the exponential morphism. 

\begin{Def}[csform][forms of complex supergroups in \emph{cs} manifolds]
  Let $G_\cplxs$ be a complex Lie supergroup with associated complex supergroup pair $(\ger g,G_{\cplxs,0})$. Let $G_0$ be a real form of $G_{\cplxs,0}$, \ie a closed subgroup whose Lie algebra is a real form of $\ger g_\ev$. The \emph{cs} Lie supergroup $G$ associated with $(\ger g,G_0)$ is called a \Define{\emph{cs} form} of $G_\cplxs$.
\end{Def}

Any \emph{cs} form $G$ of a complex Lie supergroup $G_\cplxs$ comes with a canonical morphism $G\longrightarrow G_\cplxs$ of $\cplxs$-superspaces. It can be given an expression in terms of the exponential morphism, as follows.

The real Lie supergroup $(G_\cplxs)_\reals$ associated with $G_\cplxs$ has an exponential morphism \cite{GW12}. Since its differential is complex linear for the complex structure induced by $\ger g$, it is holomorphic, and therefore induced by a unique morphism denoted by $\exp_{G_\cplxs}:\aff_{hol}(\ger g)\to G_\cplxs$. 

By the above and \thmref{Cor}{complex-prod}, the canonical morphism $G\to G_\cplxs$ is given on $T$-valued points as follows:
\[
  G\longrightarrow G_\cplxs:g\exp_G(x)\longmapsto g\exp_{G_\cplxs}(x),\quad g\in_TG_0,x\in_T\aff(\ger g_\odd)=\aff_{hol}(\ger g_\odd)
\]
Thus, the $T$-valued points of $G$ can be characterised within $G_\cplxs$ as follows. 

\begin{Prop}[csform-points]
  Let $G_\cplxs$ be a complex Lie supergroup and $G$ a \emph{cs} form of $G_\cplxs$. For any \emph{cs} manifold $T$, the $T$-valued points of $G$ are given by 
  \[
     G(T)=\Set1{g\in_TG_\cplxs}{g_0\in_{T_0}G_0}.
  \] 
\end{Prop}

\begin{proof}
  Express a $T$-valued point of $G_\cplxs$ as $g\exp_{G_\cplxs}(x)$ where $g\in_TG_{\cplxs,0}$, $x\in_T\aff_{hol}(\ger g_\odd)$. By the above considerations, $g\exp_{G_\cplxs}(x)\in G(T)$ if and only if $g\in_TG_0$.

  Thus, we need to see that
  \[
    g\in_TG_0\ \Longleftrightarrow\ g_0\in_{T_0}G_0.
  \]
  Since $(\exp_{G_\cplxs}(x))_0=1$, this will already be sufficient to prove the claim. Certainly, the left-hand statement implies that on the right-hand side. 

  On the other hand, assume that $g_0\in_{T_0}G_0$. This pins down the map underlying $g$. Local coordinates on $G_0$ can be chosen in such a way that they are the restriction of holomorphic local coordinates $(z_a)$ on $G_{\cplxs,0}$. Then $g\in_TG_0$ translates to the requirement that $g^\sharp(z_a)$ be real-valued. Since $(g^\sharp(z_a))_0=z_a\circ g_0$, this is immediate by the assumption, thereby proving the claim. 
\end{proof}

\subsection{Integration on relative \emph{cs} manifolds}\label{subs:rel-int}

Below, in \thmref{Cor}{spherical_inv}, we prove the Weyl symmetry of the spherical superfunctions $\phi_\lambda$ by the use of integral geometry. Without this fact, we would not be able to prove convergence in the delicate $c$-function expansion, nor could we derive the Harish-Chandra series, so the result is absolutely essential. 

In the proof of \loccit, we will need to handle integrals with parameters. A suitable formalism is that of fibre integrals over relative \emph{cs} manifolds. We very briefly collect the basic definitions and facts to it set up. 

\begin{Def}[rel-ber][relative Berezinian densities]
  Assume given a \emph{cs} manifold $X/Y$ over $Y$. For any fibre coordinate neighbourhood $U\subseteq X_0$, we let 
  \[
    \sh Ber_{X/Y}|_U\defi\Ber1{\sh T_{X/Y}^*(U)}.
  \]
  (Compare Ref.~\cite{manin} for the definition of the Berezinian module of a free module.) This defines a locally free $\sh O_X$-module $\sh Ber_{X/Y}$ with local basis of sections 
  \[
    Dx=D(u,\xi)=du_1\dotsm du_p\frac\partial{\partial\xi_1}\dotsm\frac\partial{\partial\xi_q}
  \]
  of parity $\equiv q\,(2)$, for any local fibre coordinate system $x=(x_a)=(u,\xi)$. Twisting by the relative orientation sheaf, we obtain $\Abs0{\sh Ber}_{X/Y}\defi or_{X_0/Y_0}\otimes_\ints\sh Ber_{X/Y}$, with corresponding local basis of sections $\Abs0{Dx}=\Abs0{D(u,\xi)}$. The local sections of the latter sheaf are called \Define{relative Berezinian densities}.
\end{Def}

\begin{Def}[rel-int][relative Berezin integral]
  Let $X/Y$ be a \emph{cs} manifold over $Y$. A retraction $r:X\to X_0$ (\ie a left inverse of $j_{X_0}:X_0\to X$) is called a \Define{retraction over $Y$} if there is a retraction $r_Y:p_X(X)\to p_X(X)_0$ \scth $r_Y\circ p_X=p_{X_0}\circ r_X$. Here, $p_X(X)\subseteq Y$ is the open subspace of $Y$ over the open set $p_{X,0}(X_0)\subseteq Y_0$.

  Fix a retraction $r$ of $X$ over $Y$. A local system $(u,\xi)$ of fibre coordinates is called \Define{$r$-adapted} if $u=r^\sharp(u_0)$ for some fibre coordinate system $u_0$ of $X_0/Y_0$.

  Let $\omega\in\Gamma(\Abs0{\sh Ber}_{X/Y})$. Then we may define 
  \[
    \fibint[Y]X[r]\omega\in \sh O_Y(p_{X,0}(X_0)),
  \]
  the \Define{relative Berezin integral} or \Define{fibre integral} of $\omega$, as follows. 

  Let $(U_i)$ be a collection of open subspaces of $X$, \scth $\supp\omega\subseteq\bigcup_iU_{i,0}$, there are $r$-adapted fibre coordinate systems $(u^i,\xi^i)$ over $Y$ defined on $U_i$, and there are $r_Y$-adapted coordinate systems $(v^i,\eta^i)$ on $V_i\defi p_X(U_i)$, \scth $(u^i,p_X^\sharp(v^i),\xi^i,p_X^\sharp(\eta^i))$ form coordinate systems of $U_i$. Let $(\chi^i)$ be a (not necessarily compactly supported) partition of unity on $\bigcup_iU_{i,0}$, subordinate to the cover. 

  We may expand 
  \[
    \omega|_{U_{i,0}}=\Abs0{D(u^i,\xi^i)}\,f^i,\quad f^i=\sum_{I,J}\xi^Ip_X^\sharp(\eta^J)r^\sharp(g^i_{IJ}),\quad g^i_{IJ}\in\sh O_{X_0}(U_{i,0}),
  \]
  with $I$ ranging through the subsets of $\{1,\dotsc,q\}$, $p|q\defi\dim_YX$, and $J$ ranging though the subsets of $\{1,\dotsc,n\}$, $m|n\defi\dim Y$. 

  With these data, we define 
  \[
    \fibint[Y]X[r]\omega\defi\sum_{i,J}\eta^Jp_Y^\sharp\Parens3{\fibint[Y_0]{X_0}\Abs0{du^i_0}\,\chi_ig^i_{\mathbf qJ}},
  \]
  whenever the integrals and the series converge absolutely. 

  Here, $\mathbf q=\{1,\dotsc,q\}$ and $\tfibint[Y_0]{X_0}\varpi$ denotes the function $y\longmapsto\smash{\int_{p_{X,0}^{-1}(y)}\varpi|_{p_{X,0}^{-1}(y)}}$, defined on $p_{X,0}(X_0)$. 
\end{Def}

\begin{Rem}
  In \thmref{Def}{rel-int}, the retraction $r_Y$ is uniquely determined by $r$. In particular, not every retraction of $X$ is over $Y$. For example, the retraction $r$ on $\aff^{2|2}$, given by $r(s)\defi(s_1+s_2s_3s_4,s_2)$ for an $T$-valued point $s\in _T\aff^{2|2}$, has no counterpart under the submersion $\psi \colon\aff^{2|2}\longrightarrow\aff^{1|2}, (s_1,s_2,s_3,s_4)\longmapsto (s_1,s_3,s_4)$.
\end{Rem}

Up to some computations in coordinates, the following fact is no harder to prove than the absolute situation where $Y=*$, \cf Ref.~\cite{AHP12}.

\begin{Th}[berint-def]
  The relative Berezin integral is well-defined independent of all choices and depends only on the choice of a retraction. In case the integrand  $\omega\in\Gamma(\Abs0{\sh Ber}_{X/Y})$ is compactly supported in the fibres over $Y$, \ie $p_{X,0}:\supp\omega\to Y_0$ is a proper map, then the integral is independent of the retraction. 
\end{Th}

Recalling the definition of the direct image with proper supports from Ref.~\cite{bredon}, the condition in the above theorem be succintly rephrased as follows:
\[
  \omega\in\Gamma(p_{X,0!}\Abs0{\sh Ber}_{X/Y}).
\]
We will use this notation in the sequel.

\begin{Cor}[change-of-var]
  Let $\vphi:X'/Y\to X/Y$ be a morphism of \emph{cs} manifolds. Let $r$ and $r'$ be retractions of $X$ and $X'$ over $Y$, respectively. Assume that $r\circ\vphi=\vphi_0\circ r'$. Then for any $\omega\in\Gamma(\Abs0{\sh Ber}_{X/Y})$, we have  
  \[
    \fibint[Y]{X'}[r']\vphi^\sharp(\omega)=\fibint[Y]X[r]\omega
  \]
  \ie both integrals exist if only one of them does, in which case they coincide. 
\end{Cor}

\subsection{Integral localization in polar coordinates}\label{subs:polar}

It is known that under suitable symmetry assumptions, there are remarkable integral localization theorems for supermanifolds \cites{dsz,sz}. Here, we show that for the special case of the supersphere, a precise form thereof can obtained by the use of polar coordinates. 

We shall make extensive use of the formalism of $T$-valued points. Compare the remarks in Subsection \ref{subs:morphisms}; the reader may also consult Ref.~\cite{a-hc}*{Appendix B}.

Denote by $y=(y_a)=(v,\eta)$ the standard coordinates on $\aff^{p|2q}$ and the Berezin--Lebesgue density by $\Abs0{D\lambda}\defi(-2\pi )^{-q} |Dy|$. The retraction associated with $y$ is
\[
  r(x)\defi x_\ev\defi(x_1,\dotsc,x_p),\quad x\in \aff^{p|2q}(T)=\Gamma\Parens1{\sh O_{T,\ev,\reals}^p\times \sh O_{T,\odd}^q}.
\] 
One should avoid to confuse $x_\ev\in\aff^p(T)$ with the underlying morphism $x_0\in \aff^p(T_0)$.

Furthermore, $\Vert \cdot\Vert ^2\colon\aff^{p|2q}\to \aff^1$ shall be given by
\begin{equation}\label{eq:normsquare}
  \Vert x\Vert ^2\defi\sum_{i=1}^px_i^2+2\sum_{j=1}^qx_{p+2j-1}x_{p+2j}  
\end{equation}
for $x\in _T\aff^{p|2q}$. Using the positive square root, this yields
\begin{align*}
  \Vert \cdot\Vert \defi\sqrt {\phantom x}\circ \Vert \cdot\Vert ^2\colon\aff^{p|2q}_{\neq0}\to \aff^1_{>0},
\end{align*}
where $\aff^{p|2q}_{\neq0}\defi\aff^{p|2q}|_{\reals^p\backslash\{0\}}$ and similarly for the subscript ``$>0$''.

\begin{Def}[rot_inv][rotationally invariant superfunctions]
  Let $f\in \Gamma(\sh O_{S\times \aff^{p|2q}})$, where $S$ is any \emph{cs} manifold. In case $p>0$, $f$ will be called \emph{rotationally invariant over $S$} if for some $\eps>0$, there exists an $f^\circ\in \sh O_{S\times\aff^1}(S_0\times(-\eps,\infty))$ such that 
  \[
    f(s,x)=f^\circ (s,\Vert x\Vert ),\quad(s,x)\in _TS\times \aff^{p|2q}_{\neq0}.
  \]
  
  In case $p=0$, $\Vert \cdot\Vert $ is not defined, so the definition of rotational invariance has to be modified as follows: $f$ is called \emph{rotationally invariant over $S$} if there exists $g\in \Gamma(\sh O_{S\times \aff^1})$, with 
  \[
    f(s,x)=g(s,\Vert x\Vert ^2),\quad (s,x)\in _T\aff^{p|2q}.
  \]
  In this case, we define $f^\circ (s,t)\defi g(s,t^2)$.
\end{Def}

\begin{Rem}
  In the case $p>0$, the super function $f^\circ|_{S_0\times(0,\infty)}$ extends to a superfunction on $S\times \aff^1$ \scth $f(s,t)=f(s,-t)$. Such an extension is given by $f^\circ(s,t)\defi f(s,te_1)$, where $t\in _T\aff^1$ and $e_1$ is the first standard basis vector of $\reals^p$. 

  Since $f^\circ$ is even in the second component, there is an extension $g$ of $f^\circ \circ ({\id_S}\times \sqrt {\phantom x})$ to $S\times\aff^1$, and $f(s,x)=g(s,\Vert x\Vert ^2)$ for all $(s,x)\in _TS\times \aff^{p|2q}$.
\end{Rem}

Here and in what follows, when handling Berezin (fibre) integrals, we will use \thmref{Th}{berint-def} and \thmref{Cor}{change-of-var} implicitly. 

\begin{Prop}[4.1]
  Let $f\in\Gamma(\sh O_{S\times \aff^{p|2q}})$ be rotationally invariant over $S$. Then
  \begin{align*}
    \fibint[S]{S\times\aff^{p|2q}}[r]\Abs0{D\lambda(x)}\,f(s,x)
      =\begin{cases}\displaystyle
          \frac{{\pi ^{\frac{p-2q}2}(-1)^q}}{\Gamma (\frac p2)}\int_0^\infty dr\,r^{\frac p2-1}\partial _r^qf^\circ (s,\sqrt r),  & p>0,
        \\  (-\pi )^{-q}\partial _{r=0}^q f^\circ (s,\sqrt r),              & p=0,
       \end{cases}
  \end{align*}
  in the sense that the integral exists if and only if the right-hand side exists, and in this case, they are equal. 
\end{Prop}

\begin{proof}
  Consider the superfunction $g$ from above. Applying Taylor expansion yields
  \begin{align*}
    g(s,t+t')
    & \equiv \sum_{k=1}^q \frac1{k!}t'^k\partial _{t'=0}^kg(s,t+t')
      \equiv \sum_{k=1}^q \frac1{k!}t'^k\partial _t^kg(s,t)\mod (t'^{q+1})
  \end{align*}
  with $t=t'=\id_{\aff^1}$. Hence, we have 
  \begin{align*}
    f(s,y)=g(s,\Vert v\Vert ^2+\Vert \eta \Vert ^2)=\sum_{k=1}^q\frac1{k!}\partial _2^kg(s,\Vert v\Vert ^2)\Vert \eta \Vert ^{2k},
  \end{align*}
  where $y=\id_{\aff^{p|q}}$, $v=\id_{\aff^p}$, $\eta =\id_{\aff^{0|q}}$, and $(\partial _2g)(s,t)\defi\partial _tg(s,t)$.
  
  The expression $\Vert \eta \Vert ^{2k}$ contains $\eta _1\dotsm \eta _{2q}$ if and only if $k=q$.
  In this case, it equals $2^q q!\eta _1\dotsm\eta _{2q}$.
  For $p=0$, this means that the integral is $(-\pi )^{-q}\partial _{t=0}^qg(s,t)$, as claimed. Similarly, if $p>0$, then it takes the form
  \begin{align*}
    \fibint[S]{S\times\aff^{p|2q}}[r]\Abs0{D\lambda(x)}\,f(s,x)
      &=(-\pi )^{-q}\int_{\reals^p}\Abs0{dv_0}\,\partial^q_2g(s,\Vert v_0\Vert ^2).\\
    \intertext{Applying polar coordinates for $p\sge2$, we obtain for $C=2\pi^{\frac{p-2q}2}(-1)^q\Gamma (\frac p2)^{-1}$}
      &=C\int_0^\infty dr\,r^{p-1}\partial ^q_2g(s,r^2)=\frac C2\int _0^\infty dr\,r^{\frac p2-1}\partial _r^qg(s,r).
  \end{align*}
  In case $p=1$, one obtains the same result by symmetry. 
\end{proof}

We obtain the following localization formula. 

\begin{Cor}[4.2]
  Let $k\sle \min(\frac p2,q)$ and $f\in\Gamma(\sh O_{S\times\aff^{p|2q}})$ be rotationally invariant and compactly supported in the fibres over $S$. Then 
  \begin{align*}
    \fibint[S]{S\times\aff^{p|2q}}[r]\Abs0{D\lambda(x)}\,f(s,x)
      =\fibint[S]{S\times\aff^{p-2k|2q-2k}}[r]\Abs0{D\lambda(x)}\,f^\circ(s,x).
  \end{align*}
\end{Cor}

\begin{proof}
  Use \thmref{Prop}{4.1} and integration by parts for $k<\frac p2$. The fundamental theorem of calculus needs to be applied for $k=\frac p2$.
\end{proof}

Combining our results, the integral of rotationally invariant superfunctions takes the following form, which depends only on $p-2q$ and not on $p|2q$. 

\begin{Cor}[4.3]
  Let $f\in \Gamma(\sh O_{S\times \aff^{p|2q}})$ be rotationally invariant and compactly supported along the fibres over $S$. Then
  \[
    \fibint[S]{S\times\aff^{p|2q}}[r]\Abs0{D\lambda(x)}\,f(s,x)=\frac{2\pi ^{\frac{p-2q}2}}{\Gamma (\frac{p-2q}2)}\int_0^\infty dr\,r^{p-2q-1}f^\circ (s,r),\quad p-2q>0.
  \]
  If $p-2q\sle 0$, then the following two cases occur:
  \begin{align*}
      \begin{cases}\displaystyle
        (-\pi )^{\frac{p-2q}2}\partial _{r=0}^{\frac{2q-p}2} f^\circ (s,\sqrt r)=(-\pi )^{\frac{p-2q}2}\frac{(\frac{2q-p}2)!}{(2q-p)!}\partial _{r=0}^{2q-p} f^\circ (s,r),        & p-2q\sle 0\ \text{even},\\  
        \displaystyle(-\pi )^{\frac{p-1-2q}2}\int_0^\infty dr\,r^{-\frac12}\partial _r^{\frac{2q+1-p}2}f^\circ (s,\sqrt r),  & p-2q<0\ \text{odd}.
       \end{cases}
  \end{align*}
\end{Cor}

\begin{proof}
  In view of \thmref{Prop}{4.1} and \thmref{Cor}{4.2}, the only case that needs some consideration is that of $p=0$. Here, Fa\`a di Bruno's formula gives 
  \[
    \partial _{r=0}^{2q}f^\circ (s,r)
       =\partial _{r=0}^{2q}f^\circ (s,\sqrt {r^2})
        =\sum_{k_1+2k_2=2q}\frac{(2q)!}{k_1!k_2!}\partial _{t=0}^{k_1+k_2}f^\circ (s,\sqrt t)\,0^{k_1}1^{k_2}.    
  \]
  All summands except for $k_2=q$ vanish.
\end{proof}

\section{Symmetric superspaces and spherical superfunctions}\label{sec:spherical}

In this section, we introduce our main objects, the \emph{spherical superfunctions}. Before proceeding to the definition of the spherical superfunctions, we collect some ancillary facts concerning symmetric superpairs and integration formul\ae{}.

\subsection{Symmetric superpairs}\label{subs:symmsup}

We review some facts on symmetric superpairs, referring to Refs.~\cites{a-hc,ASc13} for omitted details.

\begin{Def}[cartan][symmetric superpairs and Cartan decomposition]
  A symmetric superpair is a pair $(\ger g,\theta)$, where $\ger g$ be a complex Lie superalgebra and $\theta $ an involutive automorphism of $\ger g$. The eigenspace decomposition 
  \[
    \ger g=\ger k\oplus \ger p,\quad\ger k\defi\ker(1-\theta),\quad\ger p\defi\ker(1+\theta)
  \] 
  is called \emph{Cartan decomposition} of $(\ger g,\theta)$.
  
  Let $(G,K,\theta)$ be given, where $G$ is a \emph{cs} Lie supergroup with Lie superalgebra $\ger g$, $\theta$ is an involutive automorphism of $G$, and $K$ is a closed subsupergroup, $\theta|_K={\id}_K$, and its Lie superalgebra is $\ger k=\ker(1-\theta)$ (denoting the derivative of $\theta$ by the same letter). Then $(G,K,\theta)$ is called a symmetric supertriple, and the symmetric superpair $(\ger g,\theta)$ is called the \Define{infinitesimal superpair} associated with $(G,K,\theta)$.

  Then $(G,K,\theta)$ is said to admit a \Define{global Cartan decomposition} if the morphism 
  \begin{align}\label{eq:polar}
    K\times\aff(\ger p_\reals)\longrightarrow G:(k,x)\longmapsto ke^x
  \end{align}
  is an isomorphism of \emph{cs} manifolds. Here, we write $e^x\defi\exp_G(x)$.
\end{Def}

\begin{Prop}[polar]
  A pair $(G,K,\theta)$ admits a global Cartan decomposition if and only if this is true for $(G_0,K_0,\theta_0)$.
\end{Prop}

\begin{proof}
  The morphism in \eqref{eq:polar} is a local isomorphism at $(1,z)$, for any $z\in\aff(\ger p_\reals)_0$, since its derivative is given by $\ger k\times \ger p\longrightarrow\ger g:(y,x)\longmapsto y+x$. Since the morphism is $K$-equivariant, it is everywhere a local isomorphism. By the inverse function theorem \cite{leites}, it is an isomorphism if and only if the underlying map is a bijection. 
\end{proof}

\begin{Def}[red][reductive and even type conditions]

  The notions of a \Define{reductive}, \Define{strongly reductive}, or \Define{even type} symmetric superpair $(\ger g,\theta)$ are defined in Ref.~\cite{a-hc}.

  If $(\ger g,\theta)$ is the infinitesimal superpair of a symmetric supertriple $(G,K,\theta)$, then we accordingly apply these adjectives to $(G,K,\theta)$.
\end{Def}

\begin{Def}
  Let $(\ger g,\theta )$ be a reductive symmetric superpair of even type with even Cartan subspace $\ger a$.
  Then
  \begin{equation}\label{eq:root_decomp}
    \ger g=\ger m\oplus \ger a\oplus \bigoplus_{\alpha \in \Sigma }\ger g^\alpha,
    \quad\ger m\defi\ger z_\ger k(\ger a),\quad
    \ger g^\alpha\defi\bigcap_{h\in\ger a}\ker(\ad h-\alpha(h)),
  \end{equation}
  where $\Sigma\subseteq\ger a^*\setminus0$, called the set of \Define{restricted roots}, is the finite set defined by this equation. A root $\alpha\in\Sigma$ is called \Define{even} and \Define{odd}, if, respectively, $\ger g^\alpha_\ev\neq0$ and $\ger g^\alpha_\odd\neq0$. Note that roots may simultaneously be even and odd. 
  
  Given $h\in\ger a_\reals$ \scth $\alpha(h)\in\reals\setminus0$ \fa $\alpha\in\Sigma$, the subset
  \[
    \Sigma^+\defi\Set1{\alpha\in\Sigma}{\alpha(h)>0}
  \]
  is called a \Define{positive system}. Roots $\alpha\in\Sigma^+$ for which $\frac\alpha2\notin\Sigma$ are called \Define{indivisible}. We set $\vrho\defi\frac12\sum_{\alpha\in\Sigma^+}m_\alpha \alpha$, where $m_\alpha\defi\sdim\ger g^\alpha=\dim\ger g^\alpha_\ev-\dim\ger g^\alpha_\odd$.
  
  Fixing a positive system $\Sigma^+$, Equation \eqref{eq:root_decomp} takes on the form
  \[
    \ger g=\bar{\ger n}\oplus \ger m\oplus \ger a\oplus \ger n,\quad
    \ger n\defi\bigoplus_{\alpha\in\Sigma^+}\ger g^\alpha,\quad
    \bar{\ger n}\defi\theta(\ger n)=\bigoplus_{\alpha\in\Sigma^+}\ger g^{-\alpha},
  \]
  the \Define{Bruhat decomposition} of $\ger g$. Since $\ger g^\alpha\oplus\ger g^{-\alpha}$ is $\theta$-invariant, we have 
  \[
    \ger g=\ger k\oplus \ger a\oplus \ger n,
  \]
  the \Define{Iwasawa decomposition} of $\ger g$.
\end{Def}

For the remainder of this subsection, let $(G,K,\theta)$ be a reductive symmetric supertriple of even type. Fix an even Cartan subspace and a positive system of roots. Denote by $N$, $\bar N$, and $A$ the analytic subsupergroups of $G$ corresponding to $\ger n$, $\bar{\ger n}$, and $\ger a$, respectively. Here, the \Define{analytic subsupergroup} of $G$ corresponding to a subalgebra $\ger h\subseteq\ger g$ \scth $\ger h_{\ev,\reals}\defi\ger h\cap\ger g_{\ev,\reals}$ is a real form of $\ger h_\ev$ is defined to be the \emph{cs} Lie supergroup associated with the supergroup pair $(\ger h,H_0)$, where $H_0$ is the analytic (\ie, connected) subgroup of $G_0$ with the Lie algebra $\ger h_{\ev,\reals}$.

Moreover, let $M$ be the closed subsupergroup associated with the supergroup pair $(\ger m,M_0)$, where $M_0\defi Z_{K_0}(\ger a)$.

The following proposition follows in the same way as \thmref{Prop}{polar}.

\begin{Prop}[bruhat][global Bruhat decomposition]
  The morphism
  \begin{align*}
    \bar N\times M\times A\times N\longmapsto G:(\bar n,m,a,n)\longmapsto \bar nman
  \end{align*}
  is an open embedding if and only if this is true for the underlying map. In this case, we say that $(G,K,\theta)$ admits a \Define{global Bruhat decomposition}.
\end{Prop}

For the following proposition, consult Refs.~\cites{a-hc,ASc13}.

\begin{Prop}[iwasawa][global Iwasawa decomposition]
  Both of the morphisms
  \begin{align*}
    & K\times A\times N\longrightarrow G:(k,a,n)\longmapsto kan,
    \\& N\times A\times K\longrightarrow G:(n,a,k)\longmapsto nak
  \end{align*}
  are isomorphisms if and only if this is already true for one of the underlying maps. In this case, we say that $(G,K,\theta)$ admits a \Define{global Iwasawa decomposition}.
\end{Prop}

\begin{Rem}[m]
  Observe that $M$ centralises $\ger a$, since 
  \[
    \Ad(me^x)(h)=\Ad(m)\Parens1{e^{\ad(x)}(h)}=h
  \]
  for $m\in_SM_0$, $x\in_S\aff(\ger m_\odd)$, $h\in_S\aff(\ger a_\reals)$. Similarly, $M$ normalises $N$ and $\bar N$.
\end{Rem}

Let $(G,K,\theta)$ admit a global Iwasawa decomposition. We define morphisms 
\[
  k,u\colon G\to K,\quad A,H\colon G\to\aff(\ger a_\reals),\quad n,n_1\colon G\to N
\]
by requiring, for $g\in_SG$, that
\begin{equation}\label{eq:iwasawa_maps}
  g=k(g)e^{H(g)}n(g)=n_1(g)e^{A(g)}u(g).
\end{equation}
Then 
\begin{equation}\label{eq:kan-nak}
  n_1(g)=n(g^{-1})^{-1},\quad A(g)=-H(g^{-1}),\quad u(g)=k(g^{-1})^{-1}.
\end{equation}
In view of \thmref{Rem}{m}, $H$ and $k$ are right $M$-invariant, \ie
\begin{equation}\label{eq:3.1}
    H(gm)=H(g),\quad k(gm)=k(g),\quad g\in_SG, m\in_SM.
\end{equation}
Moreover, we note for later reference that 
\begin{equation}\label{eq:3.2}
  H(gh)=H\big(gk(h))+H(h),\quad k(gh)=k\Parens1{gk(h)},
\end{equation}
\fa $g,h\in_SG$. This follows from a straightforward computation on points.

The importance of the above decompositions is that they give rise to natural coordinate systems on certain homogeneous superspaces. Here and in the sequel, we will use quotients of \emph{cs} Lie supergroups. They are defined in the same way as for real Lie supergroups, see Refs.~\cites{ah-berezin,ccf} for the latter. We will need the fact that for a closed \emph{cs} Lie subsupergroup $H$ of $G$, the sheaf of superfunctions on $G/H$ is the direct image under the canonical projection $G_0\longrightarrow G_0/H_0$ of the sheaf of $H$-invariant superfunctions on $G$.

\begin{Prop}[gk-coord]
  Let $(G,K,\theta)$ admit a global Iwasawa decomposition. The morphism 
  \[
    G\longrightarrow N\times A:g\mapsto\Parens1{n_1(g),e^{A(g)}}
  \]
  induces an isomorphism $G/K\longrightarrow N\times A$. 
\end{Prop}

Let $Q=MAN$ be the closed subsupergroup of $G$ generated by $M$, $A$, and $N$.

\begin{Prop}[emb_nkm]
  Let $(G,K,\theta)$ admit global Bruhat and Iwasawa decompositions. The composite
  \begin{center}
    \begin{tikzcd}
      G\rar{k}&K\rar{}&K/M
    \end{tikzcd}
  \end{center}
  induces an isomorphism $G/Q\to K/M$. The composite 
  \begin{center}
    \begin{tikzcd}
      \bar N\rar{}&G\rar{}&G/Q\rar{\dot k}&K/M
    \end{tikzcd}
  \end{center}
  is an open embedding, which is also denoted by $k$.
\end{Prop}

\subsection{Integral formul\ae{} for supergroups and symmetric superspaces}\label{subs:intfmla}

In this subsection, we derive integral formul\ae{} for the decompositions given above. These are crucial for the proof of the central \thmref{Cor}{spherical_inv} and therefore for the proof of the $c$-function asymptotics and the Harish-Chandra series. Moreover, the proof of \loccit{} hinges on the parameter-dependent version of these formul\ae{}, which are therefore technically not a simple-minded copy of their classical relatives. We begin with some generalities, \cf Ref.~\cite{ah-berezin}.

\begin{Def}[invber][invariant Berezinian densities]
  Let $G$ be a \emph{cs} Lie supergroup and $a\colon G\times X\to X$ an action. Consider the \emph{cs} manifold $X_G\defi G\times X$ over $G$. A Berezinian density $\omega\in\Gamma(\Abs0{\sh Ber}_X)$ is said to be \Define{$G$-invariant} if
  \[
    (a,p_X)^\sharp\Parens1{p_X^\sharp(\omega)}=p_X^\sharp(\omega).
  \]
  Here, $p_X:X_G\to X$ is the projection. 
\end{Def}

The following is straightforward. 

\begin{Lem}[ber-invar]
  Let $\omega\in\Gamma(\Abs0{\sh Ber}_X)$ be a Berezinian density. Then $\omega$ is $G$-invariant if and only if, the following is true:
  \[
    \int_X \omega (x)f(gx)=\int _X \omega (x)f(x)
  \]
  for any $g\in_SG$, any \emph{cs} manifold $S,$ and any $f\in \Gamma_c(p_{X,0!}\sh O_{X_G})$.
\end{Lem}

The $G$-superspace $X$ is called \Define{analytically unimodular} if there is a non-zero $G$-invariant Berezinian density. If the action is transitive and $X$ is analytically unimodular, then $\Abs0{\sh Ber}_X$ has a global module basis, given by the choice of such a Berezinian density; moreover, it is unique up to a multiplicative constant \cite{ah-berezin}. Sufficient conditions for the analytic unimodularity of homogeneous $G$-superspaces are stated in \cite{a-hc}*{Proposition~A.2}.

\begin{Prop}[pullback_nilpot]
  Let $G$ be analytically unimodular and $\Abs0{Dg}$ a non-zero $G$-invariant Berezinian density. Let $U=-U\subseteq\aff(\ger g_{\reals})_0$ be an open neighbourhood of $0$ \scth $\exp_G:\aff(\ger g_\reals)|_U\to G$ is an open embedding. Then
  \begin{align*}
    \exp_G^\sharp(\Abs0{Dg})=\Abs0{D\lambda}(x)\,\ABer2{\frac{1-e^{-\ad x}}{\ad x}}
  \end{align*}
  on $\aff(\ger g_\reals)|_U$. Here, $\Abs0{D\lambda}$ is an adequately normalized Berezin--Lebesgue density. 
\end{Prop}

When $G$ is nilpotent, then $\str\ad x=0$ for the generic point $x$ of $\aff(\ger g_\reals)$. Hence, by \cite{a-hc}*{Proposition~A.2} and \thmref{Prop}{pullback_nilpot}, the following is immediate.

\begin{Cor}[nilpotent-haar]
  Let $G$ be a nilpotent \emph{cs} Lie supergroup. Then $G$ is analytically unimodular and 
  \[
    \exp_G^\sharp(\Abs0{Dg})=\Abs0{D\lambda}
  \]
  for an appropriate normalization of the invariant density $\Abs0{Dg}$ and the Berezin--Lebesgue density $\Abs0{D\lambda}$.  
\end{Cor}

\begin{proof}[\prfof{Prop}{pullback_nilpot}]
  The claim follows along the lines of the classical result \cite{Hel84}*{Chapter I, \S~2, Theorem 1.14}.   

  If $V$ is a sufficiently small neighbourhood of $0\in\ger g_0$, then the Campbell--Hausdorff morphism
  \[
    C:\aff(\ger g_\reals)|_U\times\aff(\ger g_\reals)|_V\to\aff(\ger g_\reals),\quad\exp_GC(x,y)=(\exp_Gx)(\exp_Gy)
  \]
  is well-defined. We let 
  \[
    X\defi\aff(\ger g_\reals)|_U\times\aff(\ger g_\reals)|_V,\quad S\defi\aff(\ger g_\reals),
  \]
  and consider $X$ as a \emph{cs} manifold over $S$ \via $p_1$. 

  Set $\vphi\defi(p_1,C)$ and let $\exp^\sharp(\Abs0{Dg})=\Abs0{D\lambda}\rho$ \fs $\rho\in\sh O_{\aff(\ger g_\reals)}(U)$. For any $f\in\Gamma_c(\sh O_{\aff(\ger g_\reals)}|_U)$, we have 
  \begin{align*}
    \int_{\aff(\ger g_\reals)}\Abs0{D\lambda}\,\rho f&=\int_G\Abs0{Dg}\,f(\log g)=\int_G\Abs0{Dg}\,f\Parens1{\log(\exp_G(-x)g)}\\
    &=\int_{\aff(\ger g_\reals)}\Abs0{D\lambda(y)}\,\rho(y)f\Parens1{C(-x,y)}\\
    &=\int_{\aff(\ger g_\reals)}\Abs0{D\lambda (y)}\,\ABer0{\sh T_{\vphi/S}} \rho \Parens0{C(x,y)}f(y).
  \end{align*}

  Here, $x$ is the generic point of $T\defi\aff(\ger g_\reals)|_U$ and the identity $y=C(x,C(-x,y))$ was applied. 
  Since $f$ was chosen arbitrarily, this implies 
  \[
    \rho (y)=\ABer1{\sh T_{\vphi/S}}(x,y)\rho\Parens0{C(x,y)},\quad y\in _T\ger g|_V.
  \]
  Setting $y=0$ shows
  \begin{align*}
    \rho(x)=\ABer1{\sh T_{\vphi/S}}(x,0)^{-1}\rho(0). 
  \end{align*}
  Hence, the claim will follow after suitable normalization and computing $\sh T_{\vphi/S}|_{y=0}$. For the latter, it suffices to prove
  \[
    T_{(x,0)/S}(\vphi)=\frac{\ad x}{1-e^{-\ad x}}
  \]
  where we consider $\ad x$ as an $\sh O_T$-linear endomorphism of $\sh O_T\otimes\ger g$. Interpreting vector fields as morphisms \via dual numbers $\mathbb D_\ev$ (see Subsection \ref{subs:tanfun}) this reduces to an equation between morphisms, which can be checked on $\aff^{0|q}$-valued points, and hence follows from the classical situation \cite{hn}.
\end{proof}

For the remainder of this subsection, let $(G,K,\theta)$ be a reductive symmetric supertriple of even type, such that a global Iwasawa decomposition exists. Further, assume that $K_0$ is compact. Then the $G\times G$-superspace $G$, the $G$-superspace $G/K$, and the $K$-superspace $K/M$ are all analytically unimodular. 

We choose corresponding invariant Berezinian densities $\Abs0{Dg}$, $\Abs0{D\dot g}$, and $\Abs0{D\dot k}$, as well as invariant densities $\Abs0{Dn}$ and $\Abs0{D\bar n}=\theta^\sharp(\Abs0{Dn})$ on $N$ and $\bar N$, respectively, and a Haar density $da$ on $A$. The normalization of these relative to each other is fixed by the Propositions \ref{Prop:int-kan}, \ref{Prop:int_k}, and \ref{Prop:pullback_barn}, which are simple generalizations of integral formul\ae{} from Refs.~\cites{a-hc,ASc13}.

We will use without further mention the fact that for compactly supported integrands, the Berezin (fibre) integral is independent of the choice of retraction (\thmref{Th}{berint-def}) and thus admits coordinate transformations (\thmref{Cor}{change-of-var}).

\begin{Prop}[int-kan][\protect{\cite{a-hc}*{Proposition 2.2}}]
  The pullback of the invariant Berezin density $\Abs0{D\dot g} $ \via the Iwasawa isomorphism is $\Abs0{D k} \otimes da\otimes \Abs0{D n} $.
  In particular,
  \begin{align*}
    \int _G\Abs0{Dg}\,f(s,g)=\int _K\Abs0{D k} \int _Ada\int _N\Abs0{D n}\,f(s,kan)e^{2\varrho (\log a)}
  \end{align*}
  for $f\in \Gamma_c(p_{1!}\sh O_{S\times G})$.
\end{Prop}

\begin{Prop}[int_k][\protect{\cite{ASc13}*{Lemma~4.2}}]
  Let $f,h\in \Gamma_c(p_{1!}\sh O_{S\times K/M})$. Then
  \begin{align*}
    \int _K\Abs0{Dk}\,f\Parens0{s,k(g^{-1}k)}h(s,k)=\int _K\Abs0{Dk}\,f(s,k)h\Parens0{s,k(gk)}e^{-2\varrho (H(gk))}
  \end{align*}
  for any $s\in_TS$ and all \emph{cs} manifolds $T$.
\end{Prop}

\begin{Cor}[int_km]
  Let $f\in \Gamma (\sh O_{S\times K/M})$. Then for any $s\in_TS$, we have 
  \begin{align*}
    \int _{K/M}\Abs0{D\dot k}\,f\big(s,k(g^{-1}k)\big)
      =\int _{K/M}\Abs0{D\dot k}\,f(s,k)e^{-2\varrho (H(gk))}.
  \end{align*}
\end{Cor}

\begin{proof}
  In the classical case, this equation follows directly from \thmref{Prop}{int_k}. In the super setting, the volume of $K$ may vanish, so one has to argue with greater care. Let $\chi\in \Gamma_c(\sh O_K)$, such that $\int_K\Abs0{Dk'}\,\chi(k')=1$. Then
  \begin{align*}
    \int_{K/M}\Abs0{D\dot k}\,f(s,k(g^{-1}k))
      &=\int_K\Abs0{Dk'}\,\chi(k')\int_{K/M}\Abs0{D\dot k}\,f\Parens1{s,k(g^{-1}k'k)}\\
      &=\int_{K/M}\Abs0{D\dot k} \int _K\Abs0{Dk'}\,\chi(k'k^{-1})f\Parens1{s,k(g^{-1}k')},
    \intertext{by the left-invariance of $\Abs0{D\dot k}$ and the right-invariance of $\Abs0{Dk'}$. By \thmref{Prop}{int_k}, this equals}
      &=\int_{K/M}\Abs0{D\dot k}\int_K\Abs0{Dk'}\,\chi\Parens1{k(gk')k^{-1}}f(s,k')e^{-2\varrho (H(gk'))},\\
    \intertext{so on applying the right-invariance of $\Abs0{Dk'}$ and the left-invariance $\Abs0{D\dot k}$ again, we obtain the equalities}
      &=\int_K\Abs0{Dk'}\int_{K/M}\Abs0{D\dot k}\,\chi\Parens1{k(gk'k)k^{-1}}f(s,k'k)e^{-2\varrho (H(gk'k))}\\
      &=\int_K\Abs0{Dk'}\int_{K/M}\Abs0{D\dot k}\,\chi\Parens1{k(gk)k^{-1}k'}f(s,k)e^{-2\varrho (H(gk))}.\\
    \intertext{Finally, the left-invariance of $\Abs0{Dk'}$ gives}
      &=\int_{K/M}\Abs0{D\dot k}\Parens3{\int_K\Abs0{Dk'}\,\chi(k')}f(s,k)e^{-2\varrho (H(gk))}\\
      &=\int_{K/M}\Abs0{D\dot k}\,f(s,k)e^{-2\varrho (H(gk))},
  \end{align*}
  which is the desired statement. 
\end{proof}

We have the following fact from Ref.~\cite{ASc13}.

\begin{Prop}[pullback_barn][\protect{\cite{ASc13}*{Proposition~4.4}}]
  The pullback of the invariant Berezin density $\Abs0{D\dot k} $ on $K/M$ \via the open embedding $k$ from \thmref{Prop}{emb_nkm} is $\Abs0{D\bar n} e^{-2\varrho (H(\bar n))}$.
\end{Prop}

\begin{Rem}[pullback_barn]
  Note that the standard retraction on $\bar N$ is in general not compatible \via $k$ with any globally defined retraction on $K/M$. Hence, \thmref{Th}{berint-def} does not apply, and the equality 
  \begin{align*}
    \int_{K/M}\Abs0{D\dot k}\,f(s,k)=\int_{\bar N}^r\Abs0{D\bar n}\,f\Parens1{s,k(\bar n)}e^{-2\varrho (H(\bar n))}
  \end{align*}
  will in general not hold for superfunctions $f\in\Gamma(\sh O_{S\times K/M})$, unless the intersection $\supp f$ with the image of $S\times\bar N$ is not compact along the fibres of $p_1$. Below, where we derive the leading asymptotics of the spherical superfunctions on $G/K$, we will resolve these issues by introducing an atlas of similar charts and cutting off the integrand in these charts by the choice of a partition of unity.

  A similar problem occurs if we use polar coordinates on $G/K$, in fact, it is more severe, leading to singularities at the boundary of the Weyl chamber. In a forthcoming paper, where we treat the inversion formula for the spherical Fourier transform in rank one, we will discuss these singularities and the `boundary terms' that they introduce at length. 
\end{Rem}

\subsection{Definition of the spherical superfunctions}\label{subs:sphdef}

Keeping the above assumptions, we can now define, for $\lambda\in\ger a^*$, the \Define{spherical superfunction} $\phi_\lambda\in\Gamma(\sh O_X)$, where $X\defi G/K$, as follows. For $g\in_TG$, where $T$ is any \emph{cs} manifold, we let   
\begin{equation}
  \phi_\lambda(g)\defi\int_{K/M}\Abs0{D\dot k}\,e^{(\lambda-\vrho)(H(gk))}.
\end{equation}
By \thmref{Lem}{ber-invar}, this defines a superfunction on $X.$

We will derive some alternative integral expressions for $\phi_\lambda$ in this subsection---in particular, we will obtain the symmetry property stated in \thmref{Cor}{spherical_inv}, which will be essential for the Harish-Chandra series expansion of $\phi_\lambda$. The statements we give here are essentially identical to the classical case, and their derivations are parallel to those given in Refs.~\cites{Hel84,Hel94}. However, they are based on \thmref{Cor}{int_km}, whose proof was a little more subtle than classically, so we briefly give the details. 

\begin{Prop}[spherical_mult]
  For any $\lambda\in\ger a^*$ and $g,h\in_TG$, we have 
  \begin{align*}
    \phi_\lambda(hg^{-1})=\int_{K/M}\Abs0{D\dot k}\,e^{(\lambda -\varrho )(H(hk))}e^{(-\lambda -\varrho )(H(gk))}
  \end{align*}
\end{Prop}

\begin{proof}
  On applying Equation \eqref{eq:3.2} twice, we obtain
  \[
    H(hg^{-1}k)=H\Parens1{hk(g^{-1}k)}+H(g^{-1}k)=H\Parens1{hk(g^{-1}k)}-H\Parens1{gk(g^{-1}k)},
  \]
  so that 
  \[
    \phi_\lambda(hg^{-1})=\int_{K/M}\Abs0{D\dot k}\,e^{(\lambda-\varrho)\Parens0{H(hk(g^{-1}k))}}e^{-(\lambda-\varrho)\Parens0{H(gk(g^{-1}k))}}.
  \]
  Setting $S\defi G\times G$ and $f(g,h,k)\defi e^{(\lambda-\varrho)\Parens0{H(hk)}}e^{-(\lambda-\varrho)\Parens0{H(gk)}}$ in \thmref{Cor}{int_km}, we arrive at our claim.
\end{proof}

Setting $h=1$ in \thmref{Prop}{spherical_mult}, the following fact is immediate. It will be instrumental in deriving a series expansion for $\phi_\lambda$.

\begin{Cor}[spherical_inv]
  For $\lambda\in\ger a^*$ and $g\in_TG$, we have $\phi_\lambda (g^{-1})=\phi _{-\lambda }(g)$.
\end{Cor}

This allows us to derive an alternative expression for $\phi_\lambda$, as follows: Consider the isomorphism $\tilde\imath_K:M\backslash K\to K/M$ induced by the inversion morphism $i_K$ of $K$. The Berezinian density
\[
  \Abs0{D\dot k_r}\defi\tilde\imath_K^\sharp(\Abs0{D\dot k}),
\]
where $\Abs0{D\dot k}$ is the $K$-invariant Berezinian density on $K/M$, is right $K$-invariant. We have the following formula, which will be used in a subsequent paper. 

\begin{Cor}
  For any $\lambda\in\ger a^*$ and $g\in_TG$, we have
  \[
    \phi_\lambda(g)=\int _{M\backslash K}\Abs0{D\dot k_r}\,e^{(\lambda+\varrho)(A(kg))}.
  \]
  where $A(\cdot)$ is defined in Equation \eqref{eq:iwasawa_maps}.
\end{Cor}

\begin{proof}
  We have by the definition of $\Abs0{D\dot k_r}$ and Equation \eqref{eq:kan-nak} 
  \begin{align*}
    \phi_\lambda(g)&=\int_{K/M}\Abs0{D\dot k} e^{(\lambda-\varrho)(H(gk))}\\
    &=\int_{K/M}\Abs0{D\dot k}e^{(-\lambda +\varrho )(A(k^{-1}g^{-1}))} 
    =\int_{M\backslash K}\Abs0{D\dot k_r}e^{(-\lambda+\varrho)(A(kg^{-1}))},
  \end{align*}
  so the assertion follows from \thmref{Cor}{spherical_inv}.
\end{proof}

\section{The Harish-Chandra \texorpdfstring{$c$}{c}-function}\label{sec:c_func}

In this section, unless something else is stated, let $(G,K,\theta)$ be a reductive symmetric supertriple of even type admitting global Iwasawa and Bruhat decompositions. A fixed even Cartan subspace $\ger a$ and positive system $\Sigma^+$ are understood, and in all the cases we will consider, the \Define{rank} $\dim\ger a-\dim\ger a\cap\dim\ger z(\ger g)$ will be one. The notation introduced in Subsection \ref{subs:symmsup} will be used freely. 

We will derive the leading asymptotics of the spherical superfunctions $\phi_\lambda$. These are governed by the \emph{Harish-Chandra $c$-function}. We compute these explicitly in terms of Euler $\Gamma$ functions. While the general structure of these functions is similar to the even setting, the occurrence of negative root multiplicities leads to a shift in the location of poles and zeroes. 

By a rank reduction procedure, this result has, in Ref.~\cite{ASc13}, been extended to arbitrary rank, on the basis of our results. In general, isotropic restricted roots occur, leading to terms in the general $c$-function formula that are not simply shifted versions of the `purely even' $c$-function.

Moreover, contrary to the `classical' setting of Riemannian symmetric spaces of the non-compact type, where the proof of convergence of the leading asymptotic term of $\phi_\lambda$ is a simple argument based on the dominated convergence theorem, in the super setting, this becomes a quite subtle issue. We need to cover the `maximal boundary' $B=K/M$ by a whole atlas of Weyl-group related `stereographical charts'. In the situations we consider (whose common feature is that they are of rank one), $B$ is a `geodesic supersphere at infinity', and this makes for two chart domains. These are essential, since the integrals on the individual domains are both divergent unless the integrand is cut off appropriately. 

For the divergences in this expansion to cancel, the Weyl symmetry of the spherical functions (\thmref{Cor}{spherical_inv}), derived above without any restriction on the rank, by the use of relative \emph{cs} manifolds and fibre integrals, will be of central importance.

Based on the derivation of the leading asymptotics, we prove the existence of a full Harish-Chandra series expansion, similar to the even case. 

We will show the existence of the $c$-function in three cases of rank one and compute it explicitly. Contrary to the classical situation, the proof of existence is by far the more difficult of these steps. Both (the proof of existence and the explicit determination) will be performed on a case-by-case basis. 

\subsection{Definition and statement of the main result}\label{subs:gk-statement}

\begin{Def}
  Let $h_0\in\ger a$ \scth $\alpha(h_0)>0$ \fa $\alpha\in\Sigma^+$. The \Define{Harish-Chandra $c$-function} is defined as follows:
  \begin{equation}\label{eq:c_func}
    c(\lambda)\defi\lim_{t\to \infty }e^{-(\lambda-\vrho)(th_0)}\phi_\lambda (e^{th_0}),
  \end{equation}
  for $\lambda \in \ger a^*$, $\Re \lambda (h_0)>0$, provided that the limit exists.

  Unless the \Define{rank} of $(G,K,\theta)$, defined to be $\dim\ger a-\dim\ger a\cap\ger z(\ger g)$, is one, it is not obvious \emph{per se} that this definition is independent of the choice of $h_0$.
\end{Def}

We will consider the cases stated in Table \ref{tab:rkone} below, where we abbreviate $\ger{gl}(p|q,\cplxs)$ by $\ger{gl}(p|q)$, \etc{}, and allow arbitrary integers $p,q\sge0$. We refer the reader to the following subsubsections for the precise definitions, in particular, regarding the \emph{cs} supergoups $G$ and $K$. 

\begin{table}[ht]\label{tab:rkone}
  \begin{center}
    \begin{tabular}{|c|c|c|c|}
      \hline
      $\ger g$&$\ger k$&$G$&$K$\\
      \hline
      \hline
      $\ger{gl}(2+p|q)$&$\ger{gl}(1)\times\ger{gl}(1+p|q)$&
      $\mathstrut\mathrm U_{cs}(1,1+p|q)$ & $\mathstrut\mathrm U(1)\times\mathrm U_{cs}(1+p|q)$\\
      \hline
      $\ger{osp}(2+p|2q)$&$\ger{osp}(1+p|2q)$&
      $\mathstrut\SOSp_{cs}^+(1,1+p|2q)$ & $\mathstrut\SOSp_{cs}(1+p|2q)$\\
      \hline
      $\ger{gl}(1|1)\times\ger{gl}(1|1)$&$\ger{gl}(1|1)$&
      $\mathstrut(\GL\times\GL)_{\cplxs|\reals}^+(1|1)$ & $\mathstrut\UGL_\reals^+(1|1)$\\
      \hline
    \end{tabular}    
  \end{center}
  \smallskip

  \caption{Riemannian symmetric superspaces under consideration}
\end{table}

For any of these symmetric pairs, we have the following statement. 

\begin{Th}[gk-fmla]
  Let be $(G,K)$ be one of the symmetric pairs listed in Table \ref{tab:rkone} and let $\alpha$ be a choice of positive indivisible restricted root. The $c$-function $c(\lambda)$ exists for $\Re\lambda>0$. For some $c_0\equiv c_0(\vrho)\neq0$, it is given by
  \begin{align*}
    c(\lambda )=
    c_0\frac{2^{-\lambda }\Gamma (\lambda )}{\Gamma \left(\frac12\Parens1{\lambda +\frac{m_\alpha}2+1}\right)\Gamma \left(\frac12\Parens1{\lambda +\frac{m_\alpha}2+m_{2\alpha }}\right)}
  \end{align*}
  if $\alpha$ is anisotropic, and if $\alpha$ is isotropic, then it is given by 
  \[
    c(\lambda)=c_0\lambda.
  \]
  Here, we identify $\lambda\equiv\frac{\Dual0\lambda\alpha}{\Dual0\alpha\alpha}$ if $\alpha$ is anisotropic and $\lambda\equiv\Dual0\lambda\alpha$ otherwise. 
\end{Th}

In particular, up to the constant $c_0$, the value of the $c$-function $c(\lambda)$ is independent of the choice of $h_0$ with $\alpha(h_0)>0$. After suitable normalisation of $h_0$, the constants $c_0$ depend on $p$ and $q$ only \via $\vrho$, see below for details. 

\subsection{Proof of Theorem \ref{Th:gk-fmla}}\label{subs:gk}

We now begin with the examination of the $c$-function for the Riemannian symmetric superspaces $X=G/K$ from Table \ref{tab:rkone}, thereby proving \thmref{Th}{gk-fmla} case-by-case. For the three individual cases, the formula is proved below in \thmref{Th}{4.5}, \thmref{Th}{c-func_sosp}, and \thmref{Cor}{c-gl11}, respectively. 

We will constantly be using the functor of points, employing \thmref{Prop}{mor-aff}, \thmref{Cor}{mor-aff-coordfree} and \thmref{Prop}{mor-aff-hol} to compute it explicitly.

\subsubsection{The unitary case}\label{ssubs:c-unitary}

In the following, let $\ger g\defi\ger{gl}(2+p|q,\cplxs)$. Moreover, we let $G\defi\mathrm U_{cs}(1,1+p|q)$ denote the \emph{cs} form of $\GL(2+p|q,\cplxs)$ corresponding to 
\[
  G_0\defi\mathrm U(1,1+p)\times \mathrm U(q),
\]
which is a the real form of the complex Lie group $G_{\cplxs,0}\defi\GL(2+p,\cplxs)\times\GL(q,\cplxs)$. 

An involution $\theta$ of $G$ resp.~of $\ger g$ is given by $\theta(x)\defi\sigma x\sigma$ where
\begin{align}\label{eq:sigma}
  \sigma \defi
    \left(\begin{array}{c:c|c}
        -\mathds{1}_1 &0      &0
      \\\hdashline
        0   &\mathds{1}_{1+p} &0
      \\\hline
        0   &0      &\mathds{1}_q
    \end{array}\right).
\end{align}

Under this involution, $\ger g$ decomposes as $\ger g=\ger k\oplus \ger p$ with
\begin{align*}
  \ger k=\left(\begin{array}{c:c|c}
              * & 0 & 0
            \\\hdashline
              0 & * & *
            \\\hline
              0 & * & *
          \end{array}\right),
  \quad
  \ger p=\left(\begin{array}{c:c|c}
              0 & * & *
            \\\hdashline
              * & 0 & 0
            \\\hline
              * & 0 & 0
          \end{array}\right).
\end{align*}
Here, dashed lines in the matrices indicate the action of the involution, whereas full lines signify the grading.

A non-degenerate invariant supersymmetric even bilinear form $b$ on $\ger g$ is given by $b(x,y)\defi\str(xy)$. We let $\ger a_\reals\subseteq \ger p_0$ be the subspace generated by the matrix
\begin{equation}\label{eq:h0-def}
  h_0\defi\left(\begin{array}{c:cc|c}
      0 & 1 & 0 & 0
    \\\hdashline
      1 & 0 & 0 & 0
    \\  0 & 0 & 0 & 0
    \\\hline
      0 & 0 & 0 & 0
  \end{array}\right).  
\end{equation}
Here and in what follows, the dimension of the rows and columns in all matrices will be $1$, $1$, $p$, and $q$, respectively. 

As ungraded vector spaces, $\ger g=\ger{gl}(2+p+q,\cplxs)$, so since $h_0$ is even, the root decomposition is the same as in the classical case. Defining $\alpha\in\ger a^*$ by $\alpha (h_0)=1$, there are exactly four roots: $\pm\alpha,\pm2\alpha$. 

The general elements of $\ger g^\alpha$ and $\ger g^{-\alpha}$ are respectively 
\begin{equation}\label{eq:galpha}
  \left(\begin{array}{c:cc|c}
            0   & 0   & B_\ev & B_\odd
          \\\hdashline
            0   & 0   & B_\ev & B_\odd
          \\  -C_\ev  & C_\ev & 0 & 0
          \\\hline
            -C_\odd  & C_\odd & 0 & 0
        \end{array}\right),\quad
  \left(\begin{array}{c:cc|c}
            0   & 0 & -B_\ev        & -B_\odd
          \\\hdashline
            0   & 0 & \phantom{-}B_\ev  & \phantom{-}B_\odd
          \\  C_\ev & C_\ev & 0         & 0
          \\\hline
            C_\odd & C_\odd & 0         & 0
        \end{array}\right).  
\end{equation}
The general elements of $\ger g^{2\alpha}$ and $\ger g^{-2\alpha}$ are respectively 
\[
  \left(\begin{array}{c:cc|c}
            -A  & A & 0 & 0
          \\\hdashline
            -A  & A & 0 & 0
          \\  0 & 0 & 0 & 0
          \\\hline
            0 & 0 & 0 & 0
        \end{array}\right),\quad
  \left(\begin{array}{c:cc|c}
            -A          & -A        & 0 & 0
          \\\hdashline
            \phantom{-}A  & \phantom{-}A  & 0 & 0
          \\  0         & 0       & 0 & 0
          \\\hline
            0         & 0       & 0 & 0
        \end{array}\right).
\]
Moreover, $\ger m$ consists of the following matrices:
\[
  \left(\begin{array}{c:cc|c}
            D & 0 & 0 & 0
          \\\hdashline
            0 & D & 0 & 0
          \\  0 & 0 & E & F
          \\\hline
            0 & 0 & G & H
        \end{array}\right).
\]

We let $\Sigma^+\defi\{\alpha,2\alpha\}$, so that $\ger n=\ger g^\alpha\oplus\ger g^{2\alpha}$ and $\bar{\ger n}=\ger g^{-\alpha}\oplus\ger g^{-2\alpha}$. From the above, we see that 
\[
  m_{\pm\alpha}=2(p-q),\quad m_{\pm2\alpha}=1,\quad\vrho=(1+p-q)\alpha.
\]
The analytic subsupergroup $K$ of $G$ with Lie superalgebra $\ger k$ is $\mathrm U(1)\times \mathrm U_{cs}(1+p|q)$. The underlying Lie group is $K_0=\mathrm U(1)\times\mathrm U(1+p)\times\mathrm U(q)$. Since the Riemannian symmetric space $G_0/K_0=\mathrm U(1,1+p)
/(\mathrm U(1)\times\mathrm U(1+p))$ is of non-compact type, $G_0$ admits global Iwasawa and Bruhat decompositions \cite{Hel62}. In view of \thmref{Prop}{iwasawa} and \thmref{Prop}{bruhat}, the same is true of $G$. 


\def\A{D}
\def\B{E}
\def\C{F}

We parametrise $\ger n$ and $\bar{\ger n}$ by setting $\tau\defi\begin{Matrix}01&1&0\\0&0&\mathds1_{p+q}\end{Matrix}^t$ and 
\begin{equation}\label{eq:XXbar}
  X(\A,\B,\C)\defi \tau
            \begin{pmatrix}
                \A & \B
              \\  \C & 0
            \end{pmatrix}
            \tau^t\sigma,\quad
  \bar X(\A,\B,\C)\defi \sigma \tau
            \begin{pmatrix}
                \A & \B
              \\  \C & 0
            \end{pmatrix}
            \tau^t,
\end{equation}
so that 
\begin{align*}
  \ger n&=\Set1{X(\A,\B,\C)}{\A\in\cplxs,\B\in\cplxs^{1\times p|q},\C\in\cplxs^{p|q\times1}},\\
  \bar{\ger n}&=\Set1{\bar X(\A,\B,\C)}{\A\in\cplxs,\B\in\cplxs^{1\times p|q},\C\in\cplxs^{p|q\times1}}.  
\end{align*}

Let $\ger g_\reals\defi\ger g_{\ev,\reals}\oplus\ger g_\odd$, $\ger g_{\ev,\reals}$ being the Lie algebra of $G_0$, and $A_\reals\defi A\cap\ger g_\reals$ for any $A\subseteq\ger g$. Then by \thmref{Prop}{csform-points}, we have an isomorphism $\vphi$, given by
\begin{gather*}\label{eq:pointsn}
    \varphi\colon\aff^1\times \aff^{p|q}\times \aff^{p|q}\longrightarrow\aff(\bar{\ger n}_\reals),\\
    (a,b,c)\longmapsto X\Parens1{\tfrac1{2i}a,\B(b,c), \C(b,c)},
\end{gather*}
on $T$-valued points $b=b_\ev+b_\odd,c=c_\ev+c_\odd\in_T\aff^{p|q}$, where we set 
\begin{gather*}
  \B(b,c)\defi
  \begin{Matrix}1
    -(b_\ev-ic_\ev)^t &
    -(b_{\odd,1}-c_{\odd,1}) & \cdots & -(b_{\odd,q}-c_{\odd,q})
  \end{Matrix}
  ,\\
  \C(b,c)\defi 
  \begin{Matrix}1
    (b_\ev+ic_\ev)^t &
    b_{\odd,1}+c_{\odd,1}& 
    \cdots &
    b_{\odd,q}+c_{\odd,q}
  \end{Matrix}^t. 
\end{gather*}
A similar statement holds for $\ger n$. 

In what follows, we will again make extensive use of the formalism of $T$-valued points, see above for explanations. For any \emph{cs} manifold $T$, the $T$-valued points of $N$ are of the form
\begin{align*}
  n_{\A\B\C}&\defi e^{X(\A,\B,\C)}=1+X(\A,\B,\C)+\tfrac12X(\A,\B,\C)^2\\
    &=\mathds1_{2+p+q}+X\Parens1{\A+\tfrac12 \B\C,\B,\C}\\
  \bar n_{\A\B\C}&\defi e^{\bar X(\A,\B,\C)}=\mathds1_{2+p+q}+\bar X\Parens1{\A+\tfrac12 \B\C,\B,\C},
\end{align*}
with $\A,\B,\C$ constrained appropriately, since $\tau^t\sigma\tau=\begin{Matrix}00&0\\0&1\end{Matrix}$. Obviously, we have the equality $\theta(n_{\A\B\C})=\bar n_{\A\B\C}$. 

\begin{Lem}[H-gl]
  The restriction of $H:G\longrightarrow\aff(\ger a_\reals)$ to $\bar N$ is given by
  \begin{align*}
    H(\bar n)
      =\frac12\log\Parens1{(1-BC)^2-4A^2}h_0,\quad
      \bar n=\bar n_{\A\B\C}\in_T\bar N.
  \end{align*}
\end{Lem}

\begin{proof}
  Formally, the calculations are the same as in the classical case, the only difference being that one has to work with generalised instead of ordinary points.

  Write $\bar n=\bar n_{\A\B\C}=ke^{th_0}n_{\A'\B'\C'}$, so that $H(\bar n)=th_0$. Then 
  \begin{equation}\label{eq:nnbar}
    n_{(-A)(-B)(-C)}\bar n_{\A\B\C}=\theta(\bar n)^{-1}\bar n=\bar n_{(-\A')(-\B')(-\C')}e^{2th_0}n_{\A'\B'\C'}.
  \end{equation}
  Let $v\defi(1,1,0)$. The simple identities 
  \[
    v\sigma\tau=0,\quad v\tau=\begin{Matrix}12&0\end{Matrix}
  \]
  give 
  \begin{align}
    vn_{\A\B\C}
      &=\begin{Matrix}11-2\A-\B\C&1+2\A+\B\C&2\B\end{Matrix},\quad v\bar n_{\A\B\C}=v\label{eq:9},\\
    n_{\A\B\C}v^t&=v^t,
    \qquad\bar n_{\A\B\C}v^t
      =\begin{Matrix}1-2\A-\B\C&1+2\A+\B\C&2\C\end{Matrix}^t\label{eq:7}.
  \end{align}
  Applying these on the right-hand side of Equation \eqref{eq:nnbar} gives
  \[
    v\theta(\bar n)^{-1}\bar nv^t=\cosh 2t+\sinh 2t=2e^{2t}.
  \]
  On the other hand, applying these on the left-hand side, we obtain
  \[
    v\theta(\bar n)^{-1}\bar nv^t=2\Parens1{(1-\B\C)^2-4\A^2}.
  \]
  This proves the lemma.   
\end{proof}

The classical approach to derive $c(\lambda)$ is to write
\begin{align*}
  \phi _\lambda (g) =\int _{K/M}\Abs0{D\dot k} e^{(\lambda -\varrho )(H(gk))}
    =\int _{\bar N}\Abs0{D\bar n} e^{(\lambda -\varrho )(H(gk(\bar n)))}.
\end{align*}
However, if $q\neq 0$, this cannot be done, since the standard retraction on $\bar N$ does not extend to $K/M$. Since $\vrho$ can be an arbitarily large negative multiple of $\alpha$, the integral on the right-hand side does not even have to exist.

The next guess might be to interchange the limit in Equation \eqref{eq:c_func} with the integral over $K/M$. However, this is not permitted, since $\exp(\lambda -\varrho )(H(e^{th_0} k e^{-th_0}))$ does admit a smooth limit, even in the classical case.

The solution we propose is to apply \thmref{Prop}{pullback_barn} after cutting off the integrand inside $k_0(\bar N_0)$. Thus, fix $\chi\in\Ct[_c^\infty]0{0,\infty}$ such that $\chi=1$ on a neighbourhood of $1$, and define $\Xi\in\Gamma(\sh O_{K/M})$ by
\begin{align*}
  \Xi(k(\bar n))\defi\chi \Parens1{(1-\B\C)^2-4\A^2},\quad\bar n=\bar n_{\A\B\C}\in_T\bar N
\end{align*}
on the open subspace $k(\bar N)\subseteq K/M$, and by zero otherwise. Clearly, this superfunction is well-defined and has compact support in the image $k(\bar N)$ of $k\colon\bar N\to K/M$. 

Then for any $w\in M_0'\defi N_{K_0}(\ger a)$, we have 
\begin{align}\label{eq:sum_int_weyl}
  \phi _\lambda (g)=\int _{K/M}\Abs0{D\dot k}\,e^{(\lambda -\varrho )(H(gk))}\Xi (k)+\int _{K/M}\Abs0{D\dot k}\,e^{(\lambda -\varrho )H(gwk)}(1-\Xi )(wk)
\end{align}
by $K$-invariance of $\Abs0{Dk}$. Denote the summands by $I_0^\lambda(g)$ and $I_\infty^\lambda(g)$, respectively.

\thmref{Prop}{pullback_barn} applies to the first of these two integrals. To do the same for the second, a good choice for $w$ has to be made. The Weyl group $W_0=M_0'/M_0=\{\pm1\}$; let $w_0\in M_0'$ be a representative of the non-trivial element.

By the Bruhat decomposition \cite{Hel62}*{Theorem~1.3}, we have 
\[
  (K_0/M_0)\setminus k_0(\bar N_0)=w_0o,\quad o\defi k(1)=1M_0.
\]
Since $1-\Xi$ vanishes in a neighbourhood of $o$, this shows that the second integrand in Equation \eqref{eq:sum_int_weyl} has compact support inside $(K_0/M_0)\setminus\{w_0o\}=k_0(\bar N_0)$.

On applying \thmref{Th}{berint-def}, \thmref{Cor}{change-of-var}, Equation \eqref{eq:3.2}, and \thmref{Prop}{pullback_barn}, we obtain
\begin{align}\label{eq:phi_barn}
  \begin{split}
    I_0^\lambda(g)&=\int _{\bar N}\Abs0{D\bar n} e^{(\lambda -\varrho )(H(g\bar n))}e^{-(\lambda +\varrho )(H(\bar n))}\Xi \Parens1{k(\bar n)}\\
    I_\infty^\lambda(g)&=\int _{\bar N}\Abs0{D\bar n} e^{(\lambda -\varrho )(H(gw_0\bar n))}e^{-(\lambda +\varrho )(H(\bar n))}(1-\Xi )\Parens1{w_0k(\bar n)}.
  \end{split}
\end{align}
For the following considerations, let $g\in_TA$. Moreover, we will identify $\ger a$ with $\cplxs$ \via $1\longmapsto h_0$ and $\ger a^*$ with $\cplxs$ \via $\lambda\longmapsto\lambda(h_0)$. It is obvious from \thmref{Cor}{exp-ad} that 
\[
  e^{th_0}\bar n_{\A\B\C}e^{-th_0}=\exp_Ge^{t\ad h_0}(\bar X(\A,\B,\C))=\bar n_{e^{-2t}\A,e^{-t}\B,e^{-t}\C},
\]
and thus, we have for $\bar n=\bar n_{\A\B\C}$:
\begin{align}
  H(e^{th_0}\bar ne^{-th_0})&=\tfrac12\log\Parens1{(1-e^{-2t}\B\C)^2-4e^{-4t}\A^2},\label{eq:han}\\
  H(e^{th_0}w_0\bar ne^{-th_0})&\begin{aligned}[t]
    &=H(e^{-th_0}\bar ne^{th_0})-2t\\
    &=\tfrac12\log\Parens1{(e^{-2t}-\B\C)^2-4\A^2},\\    
  \end{aligned}\label{eq:hawn}\\
  H(e^{th_0}k(\bar n)e^{-th_0})&=\tfrac12\log\Parens3{\frac{(1-e^{-2t}\B\C)^2-4e^{-4t}\A^2}{(1-\B\C)^2-4\A^2}},\label{eq:hak}\\
  H(e^{th_0}w_0k(\bar n)e^{-th_0})&=\tfrac12\log\Parens3{\frac{(e^{-2t}-\B\C)^2-4\A^2}{(1-\B\C)^2-4\A^2}},\label{eq:hawk}
\end{align}
where in the third and fourth line, the identity $H(ak(g)a^{-1})=H(aga^{-1})-H(g)$, which follows from Equation \eqref{eq:3.2}, was applied. 

Using Equation \eqref{eq:hak} and \thmref{Cor}{exp-ad}, one sees that 
\[
  \chi\Parens1{e^{-2H(a_tk(\bar n)a_t^{-1})}}=\chi\Parens3{\frac{(1-\B\C)^2-4\A^2}{(1-e^{-2t}\B\C)^2-4e^{-4t}\A^2}}\longrightarrow\Xi(k(\bar n))\quad(t\to\infty)
\]
for $a_t=e^{th_0}$ and $\bar n=\bar n_{\A\B\C}$. The convergence is uniform with derivatives on compact sets in $\bar n$. In particular, it holds with $k(\bar n)$ replaced by the generic point $k$ of $K$. On applying Equation \eqref{eq:hawn}, we see that 
\begin{align*}
  \Xi(w_0k(\bar n))&=\lim_{t\to\infty}\chi\Parens1{e^{-2H(a_tw_0k(\bar n)a_{\smash{t}}^{-1})}}\\
  &=\lim_{t\to\infty}\chi\Parens3{\frac{(1-\B\C)^2-4\A^2}{(e^{-2t}-\B\C)^2-4\A^2}}=\chi\Parens3{\frac{(1-\B\C)^2-4\A^2}{(\B\C)^2-4\A^2}}  
\end{align*}

We now apply the coordinates introduced by $\vphi$. We write $(s,y)$ for the coordinates and for the $T$-valued points of $\aff^1\times\aff^{2p|2q}$, where $y=(b\ c)$. In terms of the norm squared function $\Norm0\cdot^2$ defined in Equation \eqref{eq:normsquare}, we have
\[
  -4\A^2=s^2,\quad -\B\C=\Norm0y^2,\quad \bar X(\A,\B,\C)=\vphi(s,y).
\]

By \thmref{Cor}{nilpotent-haar} and since the isomorphism $\vphi$ is linear, the pullback of $\Abs0{\smash{D\dot k}}$ \via ${\exp}\circ \varphi$ is a Berezin--Lebesgue density, which, by taking $\Abs0{D\dot k}$ to be normalised adequately, may be assumed to be that from Subsection \ref{subs:polar}. We will denote the latter by $\Abs0{D\mu}$ instead of $\Abs0{D\lambda}$ to avoid confusion with the parameter $\lambda$. 

Therefore, it is immediate that Equation \eqref{eq:phi_barn} takes the form
\begin{equation}\label{eq:phiints-coord}
  \begin{split}
    I_0^\lambda(e^{th_0})e^{-t(\lambda-\varrho)}&=\int _{\aff^{1+2p|2q}}\Abs0{D\mu}\,\frac{\psi(1,e^{-t}y,e^{-2t}s)^{\frac12(\lambda -\varrho )}}{\psi(1,y,s)^{\frac12(\lambda +\varrho )}}\chi\Parens1{\psi(1,y,s)}\\
    I_\infty^\lambda(e^{th_0})e^{-t(\lambda-\varrho)}&=\int _{\aff^{1+2p|2q}}\Abs0{D\mu}\,\frac{\psi(e^{-2t},y,s)^{\frac12(\lambda -\varrho )}}{\psi(1,y,s)^{\frac12(\lambda +\varrho )}}(1-\chi )\left(\frac{\psi(1,y,s)}{\psi(0,y,s)}\right)
  \end{split}
\end{equation}
for $\psi$ defined by 
\[
  \psi(c,y,s)\defi(c+\Norm0y^2)^2+s^2.
\]
Since the integrands above are rotationally invariant over $\aff^1$, \thmref{Cor}{4.3} applies. (Here, $y$ is the variable in which we rotate, and $s$ is the base parameter.) In particular, $\phi_\lambda $ only depends on $m_\alpha$ and $m_{2\alpha}$. Therefore, if $m_\alpha>0$, we may assume that $q=0$ and use the well-known formula for the $c$-function in the classical case \cite{Hel84}*{Chapter~IV, Theorem~6.4} to arrive at the following conclusion. 

\begin{Prop}[c_u_1]
  Let $m_\alpha >0$ and $\Re \lambda >0$. Then $c(\lambda)$ exists and 
  \begin{align*}
    c(\lambda )=c_0\frac{2^{-\lambda }\Gamma (\lambda )}{\Gamma \left(\frac12\Parens1{\lambda +\frac{m_\alpha}2+1}\right)\Gamma \left(\frac12\Parens1{\lambda +\frac{m_\alpha}2+m_{2\alpha }}\right)}
  \end{align*}
  for some constant $c_0$.
\end{Prop}

In case $m_\alpha \sle 0$, we apply \thmref{Prop}{4.1}. Again, the rotation variable is $y$ and the parameter is $s$. If $f$ is one of the integrands in Equation \eqref{eq:phiints-coord}, then $f^\circ$ is determined by replacing $\psi$ with $\psi^\circ$, defined by 
\begin{equation}
   \psi^\circ(c,r,s)\defi (c+r)^2+s^2.
\end{equation}
We obtain the following intermediate result. 

\begin{Lem}[c_u_2]
  Let $m_\alpha \sle 0$ and $\Re \lambda >0$. Then 
  \begin{align}\label{eq:c_u_2}
    \begin{split}
      c(\lambda)
      &=C_0\int _0^\infty ds\,\partial _{r=0}^{-\frac{m_\alpha}2}\psi^\circ(1,r,s)^{-(\lambda +\varrho )/2}\chi \big(\psi^\circ(1,r,s)\big)
        \\&+C_0\int _0^\infty ds\,\partial _{r=0}^{-\frac{m_\alpha}2}\frac{\psi^\circ(0,r,s)^{(\lambda -\varrho )/2}}{\psi^\circ(1,r,s)^{(\lambda +\varrho )/2}}(1-\chi )\left(\frac{\psi^\circ(1,r,s)}{\psi^\circ(0,r,s)}\right)
    \end{split}
  \end{align}
  and $c(\lambda)$ is meromorphic.
\end{Lem}

\begin{proof}
The integrands have compact support, so we may interchange $\lim_{t\to \infty }$ with the integral. The same can be done for $\partial _\lambda$.
\end{proof}

In order to finally derive $c(\lambda)$, the function $\chi $ needs to be removed. The result is the following integral expression:

\begin{Prop}[4.4]
  If $m_\alpha \sle 0$, and $\Re \lambda>0$, then $c(\lambda)$ exists, and we have 
  \begin{align}\label{eq:4.4.1}
    c(\lambda )=C_0\int _0^\infty ds\,\partial _{r=0}^{1-\varrho }\big((1+r)^2+s^2\big)^{-(\lambda +\varrho)/2}.
  \end{align}
\end{Prop}

In the \emph{proof}, the following estimate will be used repeatedly.

\begin{Lem}[4.3]
  \Fa $r,s\sge 0$, $z\in \cplxs$, and $k\in\nats$, we have 
  \begin{align*}
    \Abs1{\partial_r^k\big((1+r)^2+s^2\big)^{z/2}}
      <c_k \big((1+r)^2+s^2\big)^{(\Re z-k)/2},
  \end{align*}
  where $c_k$ is some $z$-dependent constant independent of $r,s$.
\end{Lem}
\begin{proof}
  As one can see by induction,
  \begin{align*}
    \partial _r^k\big((1+r)^2+s^2\big)^{z/2}=p_k\Parens3{\frac{1+r}s}s^k \big((1+r)^2+s^2\big)^{(z-2k)/2},
  \end{align*}
  where $p_k$ is a polynomial of order at most $k$.
  Since $\lim_{t\to \infty }(t^2+1)^{-\frac k2} p_k(t)$ exists, there is a constant $c_k$ such that $|(t^2+1)^{-\frac k2}p_k(t)|<c_k$ for all $t\sge 0$. Taking $t=s^{-1}(1+r)$, this implies
  \[
    \Abs1{\partial _r^k\big((1+r)^2+s^2\big)^{z/2}}
      <c_k \big((1+r)^2+s^2\big)^{(\Re z-k)/2}.\qedhere    
  \]
\end{proof}

\begin{proof}[\prfof{Prop}{4.4}]
  Firstly, note that $-\frac{ m_\alpha}2=1-\varrho$. Denote the second integral in Equation \eqref{eq:c_u_2} by $c_{II}(\lambda)$. Since $\psi^\circ(c,ut,st)=t^2\psi^\circ(ct^{-1},u,s)$, the substitution $r=su$ in the derivative yields
  \begin{equation}\label{eq:4.3.16}
    \begin{split}
      C_0^{-1}c_{II}(\lambda)
        & =\int _0^\infty ds\,\partial ^{1-\varrho}_{u=0}s^{\varrho -1}\frac{\psi^\circ(0,su,s)^{(\lambda -\varrho)/2}}{\psi^\circ(1,su,s)^{( \lambda +\varrho)/2}}(1-\chi )\Parens3{\frac{\psi^\circ(1,su,su)}{\psi^\circ(0,su,s)}}\\
        & =\int _0^\infty ds\,\partial ^{1-\varrho }_{u=0}s^{\lambda -1}\frac{\psi^\circ(0,u,1)^{(\lambda -\varrho )/2}}{\psi^\circ(1,su,s)^{(\lambda +\varrho)/2}}(1-\chi )\Parens3{\frac{\psi^\circ(s^{-1},u,1)}{\psi^\circ(0,u,1)}}
    \end{split}
  \end{equation}
  The next step is to exchange $\int ds$ and $\partial _u^q$.
  Since only the limit $u\to 0$ is of interest, $u$ may be assumed to be small.

  Let $\varepsilon <1$ be small enough for $(1-\chi)=0$ on the interval $[1,1+4\varepsilon ^2]$.
  Since
  \[
    \frac{\psi^\circ(s^{-1},u,1)}{\psi^\circ(0,u,1)}=\frac{(s^{-1}+u)^2+1}{u^2+1}\sle 1+4\varepsilon ^2, \quad\forall s\,:\,s\eps>1,
  \]
  the function
  \begin{align*}
    [0,\varepsilon ]\times(0,\infty)\longrightarrow[0,1]:(u,s)\longmapsto (1-\chi )\left(\frac{(s^{-1}+u)^2+1}{u^2+1}\right)
  \end{align*}
  is of compact support. Therefore, $\partial ^{1-\varrho }_{u=0}$ and $\int _\delta ^\infty ds$ can be exchanged for any $\delta>0$.

  Now choose $0<\delta <1$ sufficiently small, \scth 
  \[
      (1-\chi)\Parens3{\frac{(s^{-1}+u)^2+1}{u^2+1}}=1,\quad\forall s\in(0,\delta).
  \]
  By \thmref{Lem}{4.3}, we have, for $u\sle \varepsilon $, $s\sle \delta $, and $k\sle q$, that 
  \begin{align*}
    \Abs1{s^{\lambda-1}\partial_u^k\psi^\circ(1,su,s)^{-(\lambda +\varrho)/2}}
      & < c_k s^{\Re \lambda -1+k}\psi^\circ(1,su,s)^{-(\Re \lambda +\varrho +k)/2}\\
      &\sle c_k5^{-\varrho /2}s^{\Re \lambda -1}.
  \end{align*}
  Since $s^{\Re \lambda -1}$ is integrable on $[0,\delta]$, we may also exchange $\int _0^\delta ds$ and $\partial ^{1-\varrho }_{u=0}$.

  Therefore, the right-hand side of Equation \eqref{eq:4.3.16} equals
  \begin{multline*}
    \partial ^{1-\varrho }_{u=0}
    \int _0^\infty ds\,s^{\lambda -1}\frac{\psi^\circ(0,u,1)^{(\lambda -\varrho )/2}}{\psi^\circ(1,su,s)^{(\lambda +\varrho)/2}}(1-\chi)\Parens3{\frac{\psi^\circ(s^{-1},u,1)}{\psi^\circ(0,u,1)}}\\
    =\partial ^{1-\varrho }_{u=0}\int _0^\infty \frac{ds}s
    \frac{\psi^\circ(0,\sqrt su,\sqrt s)^{(\lambda -\varrho)/2}}{\psi^\circ(1/\sqrt s,\sqrt su,\sqrt s)^{(\lambda +\varrho)/2}}(1-\chi )\Parens3{\frac{\psi^\circ(1/\sqrt s,\sqrt su,\sqrt s)}{\psi^\circ(0,\sqrt su,\sqrt s)}}.
  \end{multline*}
  We substitute $t^{-1}=\psi^\circ(0,\sqrt su,\sqrt s)=s(u^2+1)$ in the integral. Since $s^{-1}ds$ and $\psi^{\circ}(1/\sqrt s,\sqrt su,\sqrt s)=s^{-1}+2u+s(u^2+1)$ are invariant under this substitution, this leads to
  \begin{multline*}
    \partial ^{1-\varrho }_{u=0}\int_0^\infty \frac{dt}t\frac{t^{-(\lambda -\varrho )/2}}{\psi^\circ(1/\sqrt t,\sqrt tu,\sqrt t)^{(\lambda +\varrho)/2}}(1-\chi )\Parens1{t\psi^\circ(1/\sqrt t,\sqrt tu,\sqrt t)},\\
    =\partial ^{1-\varrho }_{u=0}\int _0^\infty dt\,t^{\varrho -1}\psi^\circ(1,tu,t)^{-(\lambda +\varrho)/2}(1-\chi )(\psi^\circ(1,tu,t)).
  \end{multline*}
  Again, the derivatives and the integral have to be exchanged. 
  Possibly after shrinking $\varepsilon $, one may assume that 
  \[
    (1-\chi )((1+tu)^2+t^2)=0,\quad\forall t,u:0<t,u\sle\varepsilon.
  \]
  Therefore, it suffices to consider the integral $\int _\varepsilon ^\infty dt$. Clearly, $\int _\varepsilon ^Rdt$ and $\partial ^{1-\varrho }_{u=0}$ may be exchanged, where $R$ is large enough to arrange for
  \[
     (1-\chi )((1+tu)^2+t^2)=0,\quad\forall t:t\sge R.
  \]
  On applying \thmref{Lem}{4.3} once again, we see that 
  \begin{align*}
    \Abs1{\partial _u^kt^{\varrho -1}\psi^\circ(1,tu,t)^{-(\lambda +\varrho)/2}}&<c_k t^{\varrho -1+k}\psi^\circ(1,tu,t)^{-(\Re \lambda +\varrho +k)/2}\\
      &\sle c_k5^{-\varrho /2}t^{-\Re \lambda -1},
  \end{align*}
  which is integrable over $(R,\infty)$ for $\Re \lambda >0$.

  Thus, $\int _R^\infty dt$ and $\partial ^{1-\varrho }_{u=0}$ can also be exchanged. Therefore, we find
  \begin{align*}
    c_{II}(\lambda)&=C_0\int _0^\infty dt\,\partial ^{1-\varrho }_{u=0}t^{\varrho -1}\psi^\circ(1,tu,t)^{-(\lambda +\varrho)/2}(1-\chi )(\psi^\circ(1,tu,t))\\
    &=C_0\int _0^\infty dt\,\partial ^{1-\varrho }_{r=0}\psi^\circ(1,r,t)^{-(\lambda +\varrho)/2}(1-\chi )(\psi^\circ(1,r,t))
  \end{align*}
  upon substituting $r=tu$ in the derivative. Up to a replacement of $\chi$ by $1-\chi$, this is equal to the first integral in Equation \eqref{eq:c_u_2}, and so, the assertion follows.
\end{proof}

We now arrive at our conclusion in the unitary case. 

\begin{Th}[4.5]
  The $c$-function $c(\lambda)$ for the symmetric pair of $G=\mathrm U_{cs}(1,1+p|2q)$ and $K=\mathrm U(1)\times \mathrm U_{cs}(1+p|q)$ exists for $\Re\lambda>0$. Explicitly, it is given by
  \begin{align*}
    c(\lambda )=
    c_0\frac{2^{-\lambda }\Gamma (\lambda )}{\Gamma \left(\frac12\Parens1{\lambda +\frac{m_\alpha}2+1}\right)\Gamma \left(\frac12\Parens1{\lambda +\frac{m_\alpha}2+m_{2\alpha }}\right)},\quad c_0\equiv c_0(\vrho)\neq0.
    \end{align*}
\end{Th}

\begin{proof}
  The case of $m_\alpha>0$ is the content of \thmref{Prop}{c_u_1}, so let $m_\alpha \sle 0$. Since $c(\lambda)$ is meromorphic, we may assume that $\Re \lambda>-\varrho$. \thmref{Prop}{4.4} gives 
  \[
    C_0^{-1}c(\lambda)=\int _0^\infty ds\,\partial _{r=0}^{1-\varrho }\big((1+r)^2+s^2\big)^{-(\lambda +\varrho)/2}
    =\partial ^{1-\varrho }_{r=0}\int _0^\infty ds\,\big((1+r)^2+s^2\big)^{-(\lambda +\varrho)/2}.    
  \]
  Here, integral and derivative may be exchanged due to \thmref{Lem}{4.3}, since
  \begin{align*}
    \Abs1{\partial _r^k((1+r)^2+s^2)^{-(\lambda +\varrho)/2}}
      &<c_k((1+r)^2+s^2)^{-(\Re \lambda+\varrho +k)/2}\\
      &\sle c_k (1+s^2)^{-(\Re \lambda+\varrho)/2}.
  \end{align*}
  for $k\sle 1-\varrho$. This is integrable by assumption.

  Substitution with $s=(1+r)\sqrt t$ yields
  \begin{align*}
    c(\lambda )
      &=\frac{C_0}2\partial ^{1-\varrho }_{r=0}(1+r)^{-\lambda -\varrho +1}\int _0^\infty dt\,t^{-\frac12}\big(1+t\big)^{-(\lambda +\varrho)/2}\\
      &=\frac{C_0(1-\vrho)!(-1)^{1-\varrho}}2\binom{\lambda -1}{1-\varrho }\frac{\Gamma \Parens1{\frac12}\Gamma\Parens1{\frac12(\lambda +\varrho -1)}}{\Gamma\Parens1{\frac12(\lambda+\vrho)}}\\
      & =\frac{C_0\sqrt \pi  (-1)^{1-\varrho }}2\frac{\Gamma (\lambda )\Gamma\Parens1{\frac12(\lambda +\varrho -1)}}{\Gamma (\lambda +\varrho -1)\Gamma\Parens1{\frac12(\lambda+\vrho)}}
      =C_0(-2)^{1-\varrho }\pi \frac{2^{-\lambda } \Gamma (\lambda )}{\Gamma\Parens1{\frac12(\lambda+\vrho)}^2}.
  \end{align*}
  Here, we have applied the integral formula for the Euler beta function, the identity $\binom{-x}{k}=(-1)^k\binom{x+k-1}{k}$, and the duplication formula 
  \[
    \Gamma (z)=\frac1{\sqrt \pi } 2^{z-1}\Gamma\Parens1{\tfrac12 z}\Gamma\Parens1{\tfrac12(z+1)}.  
  \]
  The claim now follows from $m_{2\alpha }=1$ and $m_\alpha =2(\varrho-1)$.
\end{proof}

\subsubsection{The ortho-symplectic case}\label{ssubs:osp}

Let $\ger g\defi\ger{osp}(1,1+p|2q,\cplxs)$ be the complex Lie subsuperalgebra of $\ger{gl}(2+p|2q,\cplxs)$, given by
\begin{align*}
  \ger{osp}(1,1+p|2q,\cplxs)=\Set1{x\in \ger{gl}(2+p|2q,\cplxs)}{x^{st^3}J+Jx=0},
\end{align*}
where 
\[
  \begin{Matrix}1R&S\\T&V\end{Matrix}^{st^3}\defi\begin{Matrix}1\phantom{-}R^t&T^t\\-S^t&V^t\end{Matrix}
\]
and
\begin{align}\label{eq:J}
  J=
    \left(\begin{array}{c:c|c}
        -\mathds 1_{1}  &0      &0
      \\\hdashline
        0   &\mathds 1_{1+p}  &0
      \\\hline
        0   &0      &J_q
    \end{array}\right).
\end{align}
Here $J_q$ denotes the $2q\times 2q$ matrix with $q$ copies of the $2\times 2$ matrix $\left(\begin{smallmatrix}0&1\\-1&0\end{smallmatrix}\right)$ on the diagonal. Therefore, $\ger g$ consists of the matrices of the form
\begin{align}\label{eq:2}
  \left(\begin{array}{c:c|c}
      0   & X_{12}  & X_{13}
    \\\hdashline
      X_{12}^t & X_{22}  & X_{23}
    \\\hline
      J_qX_{13}^t  & -J_qX_{23}^t & X_{33}
  \end{array}\right),\quad X_{22}\in \ger{so}(1+p,\cplxs),X_{33}\in \ger{sp}(2q,\cplxs),
\end{align}
with no further restrictions on the other entries. We let 
\begin{align*}
  \ger g_\reals\defi\ger g_{\ev,\reals}\oplus\ger g_\odd,\quad \ger g_{\ev,\reals}\defi\ger{so}(1+p,\reals)\times \ger{usp}(2q),
\end{align*}
where $\ger{usp}(2q)$ is the compact form of $\ger{sp}(2q,\cplxs)$, given by $\ger{usp}(2q)\defi\ger u(2q)\cap\ger{sp}(2q,\cplxs)$.

Let $G=\SOSp_{cs}^+(1,1+p|2q)$ be the \emph{cs} form of $\mathrm{OSp}(2+p|2q,\cplxs)$ given by the underlying Lie group $\SO(1,1+p,\reals)\times\mathrm{USp}(2q)$. Here, $\USp(2q)\defi\mathrm U(2q)\cap\Sp(2q,\cplxs)$.

The involution $\theta$ is the restriction of the one considered in the previous subsubsection. Let $\ger g=\ger k\oplus \ger p$ be the corresponding eigenspace decomposition. The analytic subsupergroup $K$ corresponding to $\ger k$ is $K=\SOSp_{cs}(1+p|2q)$, with underlying Lie group $\SO(1+p,\reals)\times\USp(2q)$ and Lie \emph{cs} algebra $\ger k$. For $p=0$, $K$ is the semi-direct product $\USp(2q)\ltimes\aff^{0|2q}$, where $\aff^{0|2q}$ is given the additive supergroup structure of the super-vector space $\cplxs^{0|2q}$.

Let $\ger a\subseteq \ger p_\ev$ be the even Cartan subalgebra generated by the element $h_0$ defined in Equation \eqref{eq:h0-def}. Since $\sigma $ and $J$ commute, the restricted root space decomposition for $\ger g$ is obtained by restricting that for $\ger{gl}$ from Equation \eqref{eq:galpha} to $\ger{osp}(1,1+p|2q,\cplxs)$. 

Only two restricted roots remain: $\alpha $ and $-\alpha$, where $\alpha (h_0)=1$. The general elements of $\ger g^\alpha$ and $\ger g^{-\alpha}$ are respectively 
\begin{equation}\label{eq:decomp_sosp}
      \left(\begin{array}{c:cc|c}
        0 & 0   & B_\ev & B_\odd
      \\\hdashline
        0 & 0   & B_\ev & B_\odd
      \\  B_\ev^t  & -B_\ev^t & 0 & 0
      \\\hline
        J_qB_\odd^t  & -J_qB_\odd^t  & 0 & 0
    \end{array}\right)
    ,\quad
      \left(\begin{array}{c:cc|c}
        0 & 0   & \phantom{-}B_\ev  & \phantom{-}B_\odd
      \\\hdashline
        0 & 0   & -B_\ev      & -B_\odd
      \\  B_\ev^t  & B_\ev^t    & 0     & 0
      \\\hline
        J_qB_\odd^t  & J_qB_\odd^t & 0     & 0
    \end{array}\right).
\end{equation}
The general element of $\ger m=\ger z_\ger k(\ger a)$ is
\[
      \left(\begin{array}{c:cc|c}
        0 & 0 & 0   & 0
      \\\hdashline
        0 & 0 & 0   & 0
      \\  0 & 0 & E   & F
      \\\hline
        0 & 0 & -J_qG^t  & H
    \end{array}\right).
\]
We take $\Sigma^+\defi\{\alpha\}$ as our positive system, hence $\ger n=\ger g^\alpha $ and $\bar{\ger n}=\ger g^{-\alpha }$. From Equation \eqref{eq:decomp_sosp}, we have $m_\alpha =p-2q$, so $\varrho =\frac12(p-2q)\alpha$.

Using the notation from Equation \eqref{eq:XXbar}, we set 
\[
  X(\C)\defi X(0,\B,\C),\quad\bar X(\C)\defi\bar X(0,\B,\C),\quad \B\defi\begin{Matrix}1-\C_\ev^t&-\C_\odd^tJ_q\end{Matrix}.
\]
Since the Iwasawa decomposition of $\ger g$ is the restriction of that for $\ger{gl}$ considered above, we see that 
\[
  \ger n=\Set1{X(\C)}{\C\in\cplxs^{p|2q\times1}},\quad
  \bar{\ger n}=\Set1{\bar X(\B)}{\B\in\cplxs^{p|2q\times1}}.
\]
Then by \thmref{Prop}{csform-points}, we have an isomorphism $\vphi$, given by
\[
  \vphi:\aff^{p|2q}\longrightarrow\aff(\bar{\ger n}_\reals):y\longmapsto\bar X(y),
\]
and similarly for $\ger n$. The $T$-valued points of $N$ and $\bar N$ are
\[
  n_y\defi\exp_G X(y)=\mathds1_{2+p+2q}+X(y),\quad\bar n_b\defi\mathds1_{2+p+2q}+\bar X(y),
\]
where $y\in_T\aff^{p|2q}$ and $T$ is any \emph{cs} manifold. 

As in the previous subsubsection, we will identify $\ger a$ and $\ger a^*$ with $\cplxs$ \via the bases $h_0$ and $\alpha$. The Iwasawa decompositions for $\ger{gl}$ and $\ger{osp}$ are compatible, so that we have the following lemma.

\begin{Lem}
  The restriction to $\bar N$ of the morphism $H:\bar N\to\aff(\ger a_\reals)$ is given by 
  \[
    H(\bar n)=\log(1+\Norm0y^2),\quad \bar n=\bar n_y,\quad y\in_T\aff^{p|2q},
  \]
  for any \emph{cs} manifold $T$, where $\Norm0\cdot^2$ is defined in Equation \eqref{eq:normsquare}.
\end{Lem}

\begin{proof}
  \thmref{Lem}{H-gl} applies, since $\bar n_y=\bar n_{0yb}$ where $b=-(y_\ev^t\ y_\odd^tJ_q)$. We obtain 
  \[
    H(\bar n)=\frac12\log\Parens1{(1-by)^2}=\log(1-by)=\log(1+\Norm0y^2),
  \]
  by the definition of the norm squared function. 
\end{proof}

Equations \eqref{eq:hak} and \eqref{eq:hawk} specialise to 
\begin{align*}
  H\big(e^{th_0} k(\bar n_y) e^{-th_0}\big)
    &=\log\Parens3{\frac{1+e^{-2t}\Norm0y^2}{1+\Norm0y^2}},
  \\H\big(e^{th_0} w_0k(\bar n_y) e^{-th_0}\big)
    &=\log\Parens3{\frac{e^{-2t}+\Norm0y^2}{1+\Norm0y^2}},
\end{align*}
where $w_0\in M_0'$ again denotes a representative of the non-trivial Weyl group element. 

Note that this only makes sense if $p>0$. Otherwise, the root $\alpha $ is not even and the Weyl group would be trivial in this case. Therefore, we assume $p>0$ for the time being, postponing the case of $p=0$ to the end of this subsubsection.

Let $\chi $ be as before and define $\Xi \in \Gamma (\sh O_{K/M})$ by
\begin{align*}
  \Xi \Parens1{k(\bar n_y)}\defi\chi (1+\Norm0y^2)
\end{align*}
for $\bar n_y\in _T\bar N$. As functions on $K/M$, both $\Xi $ and $k\mapsto (1-\Xi )(w_0k)$ have compact support inside $k_0(\bar N_0)$, and
\begin{align*}
  (1-\Xi)\Parens1{w_0k(\bar n_y)}=(1-\chi )\left(\frac{1+\Norm0y^2}{\Norm0y^2}\right).
\end{align*}
In analogy with Equations \eqref{eq:sum_int_weyl}, \eqref{eq:phi_barn}, and \eqref{eq:phiints-coord}, we have by \thmref{Th}{berint-def} and \thmref{Cor}{change-of-var} that 
\begin{equation}\label{eq:phi-split-osp}
  \begin{split}
    \phi _\lambda (g)&=I_0^\lambda(g)+I_\infty(g),\\
    I_0^\lambda(g)&\defi\int _{K/M}\Abs0{D\dot k}\,e^{(\lambda -\varrho )(H(gk))}\Xi (k),\\
    I_\infty^\lambda(g)&\defi\int _{K/M}\Abs0{D\dot k}\,e^{(\lambda -\varrho )H(gwk)}(1-\Xi )(w_0k),\\
    I_0^\lambda(e^{th_0})e^{-t(\lambda -\varrho)}&=\int _{\aff^{p|2q}}\Abs0{D\mu (y)} \frac{(1+e^{-2t}\Vert y\Vert ^2)^{\lambda -\varrho }}{(1+\Vert y\Vert ^2)^{\lambda +\varrho }}\chi (1+\Vert y\Vert ^2),\\
    I_\infty^\lambda(e^{th_0})e^{-t(\lambda -\varrho)}&=\int _{\aff^{p|2q}}\Abs0{D\mu (y)} \frac{(e^{-2t}+\Vert y\Vert ^2)^{\lambda -\varrho }}{(1+\Vert y\Vert ^2)^{\lambda +\varrho }}(1-\chi )\left(\frac{1+\Vert y\Vert ^2}{\Vert y\Vert ^2}\right).
  \end{split}    
\end{equation}
Here, again, $\Abs0{D\mu } $ denotes the Berezin--Lebesgue density on $\aff^{p|2q}$.

By \thmref{Cor}{4.2}, $\phi _\lambda $ and therefore $c(\lambda )$ only depends on $m_\alpha $ and not on $p$ and $q$ separately in this case, too. Therefore, for $m_\alpha >0$, \cite{Hel84}*{Chapter~IV, Theorem~6.4} implies the following. 

\begin{Prop}[c_osp_1]
  Let $m_\alpha >0$. Then $c(\lambda)$ exists for $\Re \lambda >0$, and we have 
  \begin{align*}
    c(\lambda )=c_0\frac{2^{-\lambda }\Gamma (\lambda )}{\Gamma \left(\frac12(\lambda +\frac{m_\alpha}2+1)\right)\Gamma \left(\frac12(\lambda +\frac{m_\alpha}2+m_{2\alpha })\right)},\quad c_0\equiv c_0(\vrho)\neq0.
  \end{align*}
\end{Prop}

For $2\varrho=m_\alpha\sle 0$, we get a similar expression. 

\begin{Prop}[c_osp_2]
  Let $m_\alpha \sle 0$ and $p>0$. Then $c(\lambda)$ exists for $\Re \lambda >0$, and  
  \begin{align*}
    c(\lambda )=c_0'\frac{\Gamma (\lambda )}{\Gamma (\lambda +\varrho )},\quad c_0'\equiv c_0'(\vrho).
  \end{align*}
\end{Prop}

\begin{proof}
  Following \thmref{Cor}{4.3}, we distinguish between the even and odd cases. Firstly, assume that $m_\alpha$ be even. In this case, the integrals in Equation \eqref{eq:phi-split-osp} reduce to a pointwise derivative, which may be exchanged with limits to arrive at 
  \[
      c(\lambda )
      = C\partial ^{-\varrho }_{r=0}\left((1+r)^{-(\lambda +\varrho )}\chi (1+r)+r^{\lambda -\varrho }(1+r)^{-(\lambda +\varrho )}(1-\chi )\left(\frac{1+r}r\right)\right).       
  \]    
  for some constant $C$. As functions of $r$, $\chi (1+r)$ and $(1-\chi )(r^{-1}(1+r))$ are constant near zero. Furthermore, $\partial ^k_{r=0}r^{\lambda-\varrho }=0$ for $k\sle -\varrho $, so that 
  \[
    c(\lambda )  =C q!\binom{-(\lambda +\varrho )}{-\varrho }
        =C (-1)^\varrho  q!\binom{\lambda -1}{-\varrho }
        =C (-1)^\varrho  \frac{\Gamma (\lambda )}{\Gamma (\lambda +\varrho )}.    
  \]

  Now, assume that $m_\alpha<0$ is odd. By \thmref{Lem}{4.3}, we may exchange limits and integrals, to obtain 
  \begin{equation}\label{eq:c_osp_3_1}
    \begin{split}
      c(\lambda )
        =C&\int _0^\infty dr\,r^{-\frac12}\partial _r^{\frac12-\varrho}\Parens1{(1+r)^{-(\lambda +\varrho )}\chi (1+r)}\\
        &+C\int _0^\infty dr\,r^{-\frac12}\partial _r^{\frac12-\varrho}\Parens1{r^{\lambda -\varrho }(1+r)^{-(\lambda +\varrho )}(1-\chi )\Parens1{r^{-1}(1+r)}}.
    \end{split}
  \end{equation}
  
  Denote the two integrals by $c_I(\lambda)$ and $c_{II}(\lambda)$, respectively. As in the proof of \thmref{Prop}{4.4}, we will rewrite these in order to show that the respective contributions of the cutoff function $\chi$ cancel. 

  This will be done by integrating by parts. Recall that $(1-\chi )(r^{-1}(1+r))=1$ for small $r$. Therefore, for $\Re\lambda >0$ and $k\sle \frac12-\varrho $, we have 
  \begin{align*}
    \lim_{r\to 0}&(\partial _r^{k-1}r^{-\frac12})\partial _r^{-\varrho -\frac12-k}\Parens1{r^{\lambda -\varrho }(1+r)^{-(\lambda +\varrho )}(1-\chi )\Parens1{r^{-1}(1+r)}}
      \\& =\lim_{r\to 0}r^{\frac12-k}\sum\nolimits_{\ell=0}^{\frac12-\varrho -k}c_\ell r^{\lambda -\varrho -(\frac12-\varrho -k-\ell)}(1+r)^{-(\lambda -\varrho )-\ell}
      \\&=\lim_{r\to 0}\sum\nolimits_{\ell=0}^{\frac12-\varrho -k}c_\ell r^{\lambda+\ell}(1+r)^{-(\lambda-\varrho)-\ell}
        =0
  \end{align*}
  
  Therefore, no boundary terms occur upon integrating by parts, and 
  \begin{align*}
    c_{II}(\lambda)&=\Parens1{\tfrac12-\vrho}!\binom{-\vrho}{\tfrac12-\vrho}\int _0^\infty\frac{dr}r\,r^\lambda(1+r)^{-(\lambda +\varrho )}(1-\chi )\Parens1{r^{-1}(1+r)}\\
    &=\Parens1{\tfrac12-\vrho}!\binom{-\vrho}{\tfrac12-\vrho}\int _0^\infty\frac{ds}s\,s^{-\lambda }\Parens1{s^{-1}(1+s)}^{-(\lambda +\varrho )}(1-\chi )(1+s)\\
    &=\Parens1{\tfrac12-\vrho}!\binom{-\vrho}{\tfrac12-\vrho}\int_0^\infty\frac{ds}s\, s^\varrho(1+s)^{-(\lambda +\varrho )}(1-\chi )(1+s).
  \end{align*}
  where $r^{\vrho-1}r^{\lambda-\vrho}=r^{\lambda-1}$ and the substitution $s=r^{-1}$ were applied.

  Since similarly, for $k\sle \frac12-\varrho $, we have 
  \[
    \lim_{r\to \infty }\partial _r^{k-1} r^{-\frac12}\partial _r^{\frac12-\varrho -k}\Parens1{(1+r)^{-(\lambda +\varrho )}(1-\chi )(1+r)}
    \simeq \lim_{r\to \infty }r^{\frac12-k}(1+r)^{-\lambda -\frac12+k}
    =0,
  \]
  we may integrate by parts to find 
  \begin{multline*}
    \int _0^\infty dr\,r^{-\frac12}\partial _r^{\frac12-\varrho }\Parens1{(1+r)^{-(\lambda +\varrho )}(1-\chi )(1+r)}\\
    =\Parens1{\tfrac12-\vrho}!\binom{-\vrho}{\tfrac12-\vrho}\int_0^\infty\frac{ds}s\, s^\varrho(1+s)^{-(\lambda +\varrho )}(1-\chi )(1+s),
  \end{multline*}
  since $(1-\chi )(1+r)$ vanishes for $r$ small and equals $1$ for $r$ large.
  
  Combining $c_I(\lambda)$ and $c_{II}(\lambda)$ leads to 
  \[
    c(\lambda )
       =C\int _0^\infty dr\,r^{-\frac12}\partial _r^{\frac12-\varrho }(1+r)^{-(\lambda +\varrho )}
      =C(-1)^{-\varrho -\frac12}\sqrt\pi \frac{\Gamma (\lambda )}{\Gamma (\lambda +\varrho )}
  \]
  for $\Re\lambda>-\vrho$, as in the proof of \thmref{Th}{4.5}. Moreover, Equation \eqref{eq:c_osp_3_1} shows that $c(\lambda)$ is holomorphic in $\lambda $ for $\Re\lambda>0$, since the integrands are sufficiently bounded. Therefore, the case of $0<\Re \lambda \sle -\varrho$ follows by analytic continuation.
\end{proof}

The case where $p=0$ is still open. This situation is interesting, in that there are only purely odd restricted roots (which, however, are anisotropic). Fortunately, one can evaluate $\phi_\lambda$ explicitly in this case. 

\begin{Prop}[spherical_osp_p0]
  If $p=0$, then $\phi _\lambda $ is for some constant $c_1$ given by
  \begin{align*}
    \phi _\lambda (e^{th_0})=c_1e^{\lambda t}\sum_{k=0}^{-\varrho }\binom{\lambda -\varrho }{k}\binom{-\lambda -\varrho }{-\varrho -k}e^{(-\varrho -2k)t}
  \end{align*}
\end{Prop}

\begin{proof}
  The \emph{cs} manifolds $K/M$ and $\bar N$ admit only one retraction, since the underlying spaces are trivial. Hence, we may by \thmref{Th}{berint-def} and \thmref{Cor}{change-of-var} pull back the defining integral of $\phi_\lambda$ and apply Equation \eqref{eq:3.2}, to obtain 
  \begin{align*}
    \phi _\lambda (e^{th_0})e^{-t(\lambda -\varrho )}
      & =\int _{\bar N}\Abs0{D\bar n}\,e^{(\lambda -\varrho )(H(e^{th_0}\bar ne^{-th_0}))}e^{-(\lambda +\varrho )(H(\bar n))}
      \\& =\int _{\aff^{0|2q}}\Abs0{D\mu (y)}\,(1+e^{-2t}\Vert y\Vert ^2)^{\lambda -\varrho }(1+\Vert y\Vert ^2)^{-(\lambda +\varrho )}
      \\& \simeq\partial _{r=0}^q(1+e^{-2t}r)^{\lambda -\varrho }(1+r)^{-(\lambda +\varrho )}\\
      &\simeq q!\sum\nolimits_{k=0}^q\binom{\lambda -\varrho }{k}\binom{-\lambda -\varrho }{q-k}e^{-2kt}.
  \end{align*}
  This proves the claim. 
\end{proof}

\begin{Cor}[c_osp_4]
  If $p=0$, the $c(\lambda)$ exists for $\Re\lambda>0$ and 
  \begin{align*}
    c(\lambda )=c_0'\frac{\Gamma (\lambda )}{\Gamma (\lambda +\varrho )},\quad c_0'\equiv c_0'(\vrho).
  \end{align*}
\end{Cor}

Using the duplication formula, \thmref{Prop}{c_osp_1}, \thmref{Prop}{c_osp_2}, and \thmref{Cor}{c_osp_4} combine to the following result. 

\begin{Th}[c-func_sosp]
 The $c$-function $c(\lambda)$ for the symmetric pair of the \emph{cs} Lie supergroups $G=\SOSp^+_{cs}(1,1+p|2q)$ and $K=\SOSp_{cs}(1+p|2q)$ exists for $\Re\lambda>0$. Explicitly, 
  \begin{align*}
    c(\lambda )=
    c_0\frac{2^{-\lambda }\Gamma (\lambda )}{\Gamma \left(\frac12\Parens1{\lambda +\frac{m_\alpha}2+1}\right)\Gamma \left(\frac12\Parens1{\lambda +\frac{m_\alpha}2+m_{2\alpha }}\right)},\quad c_0\equiv c_0(\vrho)\neq0.
  \end{align*}
\end{Th}

\subsubsection{The case of $\mathrm{GL}(1\vert 1)$}\label{ssubs:gl11}

Let $\ger g\defi \ger{gl}(1|1,\cplxs)\times\ger{gl}(1|1,\cplxs)$. Then $\ger g_\ev$ is Abelian. We will write the elements of $\ger g$ as double matrices of the form
\begin{align*}
  \left(\begin{array}{cc|cc}
      A & B & E & F
    \\  C & D & G & H
  \end{array}\right).
\end{align*}
Let $\ger g_\reals\defi\ger g_{\ev,\reals}\oplus\ger g_\odd$ be given by requiring $\ger g_{\ev,\reals}$ to be the Lie algebra of the Lie group $G_0$ whose general elements are 
\[
  \left(\begin{array}{cc|cc}
      z & 0 & \bar z^{-1} & 0
    \\  0 & r & 0   & s
  \end{array}\right),\quad z\in \cplxs^\times,r,s>0.
\]
The \emph{cs} form $G$ of $G_\cplxs=\GL(1|1,\cplxs)\times\GL(1|1,\cplxs)$ defined by $G_0$ will be denoted by $(\GL\times\GL)_{\cplxs|\reals}^+(1|1)$. An involution $\theta $ on $\ger g$ can be defined by
\begin{align*}
  \theta \left(\begin{array}{cc|cc}
      A & B & E & F
    \\  C & D & G & H
  \end{array}\right)
  =
  \left(\begin{array}{cc|cc}
      E & F & A & B
    \\  G & H & C & D
  \end{array}\right).
\end{align*}
The general elements of the eigenspaces $\ger k$ and $\ger p$ are respectively
\[
  \left(\begin{array}{cc|cc}
      A & B & A & B
    \\  C & D & C & D
  \end{array}\right),\quad
  \left(\begin{array}{cc|cc}
      A & B & -A  & -B
    \\  C & D & -C  & -D
  \end{array}\right).  
\]
Let $K$ be the analytic \emph{cs} Lie subgroup of $G$ with Lie \emph{cs} algebra $\ger k$, denoted by $K=\UGL_\reals^+(1|1)$. The general element of $K_0$ is 
\[
  \left(\begin{array}{cc|cc}
    z & 0 & z & 0
  \\  0 & r & 0 & r
  \end{array}\right),\quad z\in\UU(1),r>0.
\]
Note that although $K_0$ is not compact, its adjoint image in $\GL(\ger p_{\ev,\reals})$ is. This explains why there is no contradiction in the fact that the negative of the supertrace form is actually positive on $\ger p_{\ev,\reals}$, so that $(G_0,K_0)$ is a Riemannian symmetric pair. 

The only choice of Cartan subalgebra is $\ger a\defi\ger p_\ev$, since $\ger p_\ev$ is Abelian and consists of semisimple elements. Let
\begin{align*}
  h_1\defi 
      \left(\begin{array}{cc|cc}
          1 & 0 & -1  & 0
        \\  0 & 0 & 0 & 0
      \end{array}\right),\quad
  h_2\defi 
      \left(\begin{array}{cc|cc}
          0 & 0 & 0 & 0
        \\  0 & 1 & 0 & -1
      \end{array}\right)
\end{align*}
be a basis of $\ger a$. An easy calculation shows that there are only the roots $\pm\alpha$, with $\alpha$ defined by 
\[
  \alpha (a_1h_1+a_2h_2)\defi a_1-a_2.  
\]
We have $\ger m=\ger z_\ger k(\ger a)=\ger k$. The general elements of $\ger g^\alpha$ and $\ger g^{-\alpha}$ are respectively 
\[
  \left(\begin{array}{cc|cc}
      0 & B & 0 & 0
    \\  0 & 0 & G & 0
  \end{array}\right),\quad  
  \left(\begin{array}{cc|cc}
      0 & 0 & 0 & F
    \\  C & 0 & 0 & 0
  \end{array}\right).
\]

We let $\Sigma^+\defi\{\alpha\}$ be our positive system, so $\ger n=\ger g^\alpha $, $\bar{\ger n}=\ger g^{-\alpha }$, $m_\alpha =-2$, and $\varrho =-\alpha$. Since $G_0$ admits a global Iwasawa decomposition, this is also the case for $G$. Indeed, every element in $G_0$ can be decomposed uniquely as 
\[
  \left(\begin{array}{cc|cc}
      z & 0 & \bar z^{-1} & 0
    \\  0 & r & 0   & s
  \end{array}\right)
  =\left(\begin{array}{cc|cc}
      \frac{z}{|z|} & 0 & * & * \\ 
      0   & \sqrt{rs}  & * & *
  \end{array}\right)
  \left(\begin{array}{cc|cc}
      \log|z| & 0     & * & * \\
    0 & \frac12\log(\frac rs) & * & *
  \end{array}\right),
\]
where the first factor is in $K_0$ and the second in $A_0$. Note that $N_0=1$.

For any \emph{cs} manifold $T$, the $T$-valued points of $N$ and $\bar N$ are of the form
\begin{align*}
  n_{B,G}=
    \left(\begin{array}{cc|cc}
        1 & B & 1 & 0
      \\  0 & 1 & G & 1
    \end{array}\right),\quad
  \bar n_{C,F}=
    \left(\begin{array}{cc|cc}
        1 & 0 & 1 & F
      \\  C & 1 & 0 & 1
    \end{array}\right),\quad B,C,F,G\in \Gamma (\sh O_{T,\odd}).
\end{align*}

\begin{Lem}[h-gl11]
  The restriction to $\bar N$ of the morphism $H:G\to\aff(\ger a_\reals)$ is given by 
  \[
    H(\bar n)=\tfrac12 CFh^+,\quad \bar n=\bar n_{C,F}\in_T\bar N,
  \]
  where we set $h^\pm\defi h_1\pm h_2$.
\end{Lem}

\begin{proof}
  Again, $H\colon \bar N\longrightarrow \ger a$ will be derived from the equation
  \begin{align*}
    \theta (\bar n_{C,F})^{-1}\bar n_{C,F}=\theta (n_{B,G})^{-1}\exp(2H(\bar n_{C,F})) n_{B,G},\qquad \bar n_{C,F}\in _T\bar N, n_{B,G}\in_T N.
  \end{align*}
  The left-hand side equals
  \begin{align*}
        \left(\begin{array}{cc|cc}
            1 & -F  & 1 & 0
          \\  0 & 1 & -C  & 1
        \end{array}\right)
        \left(\begin{array}{cc|cc}
            1 & 0 & 1 & F
          \\  C & 1 & 0 & 1
        \end{array}\right)
      =
        \left(\begin{array}{cc|cc}
            1-FC  & -F  & 1 & F
          \\  C    & 1   & -C & 1-CF
        \end{array}\right).
  \end{align*}
  Upon writing $H(\bar n_{C,F})=h_1t_1+h_2t_2$ for some $t_1,t_2\in \Gamma (\sh O_{T,\ev})$, the right-hand side takes the form 
  \begin{align*}
        \bigg(\begin{array}{cc|}
            1 & 0\\
            -G  & 1
        \end{array}
        &
        \begin{array}{cc}
          1 & -B\\
          0 & 1
        \end{array}
        \bigg)
        \left(\begin{array}{cc|cc}
          e^{2t_1}  & 0 & e^{-2t_1} & 0
        \\  0   & e^{2t_2}  & 0     & e^{-2t_2}
        \end{array}\right)
        \left(\begin{array}{cc|cc}
          1 & B & 1 & 0
        \\  0 & 1 & G & 1
        \end{array}\right)\\
        &=
        \left(\begin{array}{cc|cc}
            e^{2t_1}    & Be^{2t_1}     & e^{-2t_1}-BGe^{-2t_2} & -Be^{-2t_2}
          \\  -Ge^{2t_1}  & e^{2t_2}-GBe^{2t_1} & Ge^{-2t_2}      & e^{-2t_2}
        \end{array}\right),
  \end{align*}
  Therefore, $e^{2t_1}=1-FC$ and $e^{-2t_2}=1-CF=1+FC$, so that
  \[
    t_1=t_2=-\tfrac12FC=\tfrac12CF,
  \]
  which proves the lemma. 
\end{proof}

We find the following formula for $\phi_\lambda$.

\begin{Prop}[phi_iso]
  For a suitable normalisation of the invariant density on $K/M$, we have for any $\lambda \in \ger a^*$
  \begin{align*}
    \phi _\lambda (e^h)=-\lambda (h^+)e^{\lambda (h)}\sinh\alpha (h),\quad h\in_T\aff(\ger a_\reals).
  \end{align*}
\end{Prop}

Remarkably, the function $\phi_\lambda$ is not invariant under the symmetry $h\longmapsto-h$. This is due to the fact that there are no even restricted roots in the case under consideration, so that the even Weyl group is trivial. 

\begin{proof}[\prfof{Prop}{phi_iso}]
  We introduce coordinates $\xi_1,\xi_2$ on $\bar N$ by
  \begin{align*}
        \xi _1(\bar n_{C,F})=C,
    \quad \xi _2(\bar n_{C,F})=F.
  \end{align*}
  In terms of these coordinates, we may take $\Abs0{D\bar n} =\Abs0{D\xi }=\Abs0{D(\xi_1,\xi_2)}$.

  The spaces $K/M$ and $\bar N$ being of purely odd dimension, they both admit only one retraction. Thus, we may by \thmref{Th}{berint-def} and \thmref{Cor}{change-of-var} pull back the defining integral of $\phi_\lambda$, leading to 
  \begin{align}\label{eq:phi_iso}
    \phi _\lambda (e^h)e^{-(\lambda -\varrho )(h)}
      =\int _{\bar N}\Abs0{D\bar n} e^{(\lambda -\rho )(H(e^h\bar ne^{-h}))}e^{-(\lambda +\varrho )(H(\bar n))}
  \end{align}
  for $h\in_T\aff(\ger a_\reals)$, by \thmref{Prop}{pullback_barn} and Equation \eqref{eq:3.2}. 

  Clearly, we have by \thmref{Cor}{exp-ad} that 
  \[
    e^h\bar n_{C,F}e^{-h}=\bar n_{e^{-2\alpha (h)}C,e^{-2\alpha (h)F}},
  \] 
  and $\varrho (h^+)=-\alpha (h^+)=0$. Furthermore, \thmref{Lem}{h-gl11} gives 
  \begin{align*}
    e^{\lambda (H(\bar n_{C,F}))}
      =e^{\frac12\lambda (h^+)CF}
      =1+\tfrac12\lambda (h^+)CF.
  \end{align*}
  Inserting these findings into Equation \eqref{eq:phi_iso}, we arrive at 
  \begin{align*}
    \phi _\lambda (e^h)e^{-(\lambda -\varrho )(h)}
    & =\int _{\bar N}\Abs0{D\xi}\,\Parens1{1+\tfrac12\lambda(h^+)e^{-2\alpha (h)}\xi _1\xi_2}\left(1-\tfrac12\lambda(h^+)\xi_1\xi_2\right)\\
    &=\tfrac12\lambda (h^+)(e^{-2\alpha (h)}-1)=-\lambda(h^+)e^{\vrho(h)}\sinh\alpha(h),
  \end{align*}
  proving the assertion.
\end{proof}

We immediately obtain the following expression for the $c$-function.

\begin{Cor}[c-gl11]
  Let $\lambda\in\ger a^*$ be arbitrary. Using $h_0=c_+h^++c_-h^-$ in Equation \eqref{eq:c_func}, $c(\lambda)$ exists if and only if $c_-=0$ or $\Re c_->0$, and then 
  \begin{align*}
    c(\lambda)=
    \begin{cases}
      0, &c_-=0,\\
      -\tfrac12\lambda (h^+), &\Re c_->0.
    \end{cases}
  \end{align*}
\end{Cor}

\begin{proof}
  Since $\alpha(h^+)=0$ and $\alpha(h^-)=2$, we have 
  \[
    e^{\vrho(th_0)}\sinh\alpha(th_0)=e^{-2c_-t}\sinh 2c_-t=\tfrac12\Parens1{1-e^{-4c_-t}}.
  \]
  Thus, the limit for $t\to\infty$ exists if and only if $c_-=0$ or $\Re c_->0$. In these cases, it equals $0$ and $\tfrac12$, respectively. 
\end{proof}

Thus, in particular, $c(\lambda)$ can be defined independently of a choice of $h_0$ with $\alpha(h_0)=2c_->0$. Since $\alpha(h)=\tfrac12\Dual0{h^+}h$, the above result for $c(\lambda)$ can be rewritten in the form stated in \thmref{Th}{gk-fmla}.

\section{Determination of the spherical functions}\label{sec:hcseries}

In this subsection, we derive a convergent series expansion for the function $\phi_\lambda$. The standard procedure to do so is to solve a differential equation by a pertubation ansatz, \cf \cite{Hel84}*{Chapter~IV, §5}. We follow the same procedure. 

In the following, suppose $(G,K)$ is one of the two following symmetric pairs
\begin{equation}\label{eq:symmpairs-exgl11}
  \begin{gathered}
    \Parens1{\UU_{cs}(1,1+p|2q),\UU(1)\times\UU_{cs}(1+p|2q)},\\
    \Parens1{\SOSp_{cs}^+(1,1+p|2q),\SOSp_{cs}(1+p|2q)}.    
  \end{gathered}
\end{equation}
For $(G,K)=\Parens1{(\GL\times\GL)_{\cplxs|\reals}^+(1|1),\UGL_\reals^+(1|1)}$, a closed expression for $\phi_\lambda$ was derived above in \thmref{Prop}{phi_iso}. Moreover, in the orthosymplectic case, a closed expression was obtained for $p=0$, in \thmref{Prop}{spherical_osp_p0}. 
This section aims to remove such restrictions.

Let $h_0\in\ger a_\reals$ be determined uniquely by $\alpha(h_0)=1$. By the use of the bases $h_0$ and $\alpha$ of $\ger a$ and $\ger a^*$, respectively, we will identify these spaces with $\cplxs$ where convenient.

\begin{Prop}[phi_laplace]
  Let $\Delta (L)$ be the differential operator on $A$ given by
  \begin{align*}
    \Delta (L)(f)(e^{th_0})\defi\Parens1{\partial _t^2+(m_\alpha \coth(t)+2m_{2\alpha}\coth(2t))\partial _t}f(e^{th_0}).
  \end{align*}
  Then, for any $\lambda\in\ger a^*$, $\phi_\lambda$ is an eigenfunction of $\Delta(L)$. More precisely, 
  \begin{align}\label{eq:phi_laplace}
    \Delta (L)\phi _\lambda =(\lambda ^2-\varrho ^2)\phi _\lambda .
  \end{align}
\end{Prop}

\begin{proof}
  One follows the standard procedure, showing that $\Delta(L)$ is the $N$-radial part of the Laplacian on $G/K$, \cf \cite{GV88}*{§4.2}. This can also be seen by a computation in $\Uenv0{\ger g}$ along the lines of the proof of \cite{warner-vol2}*{Proposition 9.1.2.11}. The eigenvalue equation is derived as in Ref.~\cite{a-hc}.
\end{proof}

Since $e^{(\lambda -\varrho)t}$ is an eigenfunction of $\partial_t(\partial_t+2\vrho)$ for the correct eigenvalue, we make the following perturbation ansatz:
\begin{align}\label{eq:def_Phi}
  \Phi _\lambda (e^{th_0})=e^{(\lambda -\varrho )t}\sum\nolimits_{\ell=0}^\infty \gamma_\ell(\lambda )e^{-2\ell t},\quad e^{th_0}\in A^+,
\end{align}
where $A^+$ is the open subset of $A_0=A$ defined by $A^+\defi\exp((0,\infty)h_0)$. For $\eps>0$, we denote by $A_\eps^+$ the subset of $A_0$ defined by $A_\eps^+\defi\exp([\eps,\infty)h_0)$.

\begin{Prop}[Phi-conv]
  Let $\lambda\in\ger a^*$, $\lambda\equiv\lambda(h_0)\notin\nats$. Assuming that the series $\Phi_\lambda$ in Equation \eqref{eq:def_Phi} converges to a solution of Equation \eqref{eq:phi_laplace}, absolutely on $A_\eps^+$ for any $\eps>0$, the coefficients $\gamma_\ell(\lambda)$ are determined by the choice of $\gamma_0(\lambda)$ and the relation 
  \begin{equation}\label{eq:coeff_Phi}
    \begin{split}
    \ell(\ell-\lambda)\gamma_\ell(\lambda)=&\;
      \frac{m_\alpha}2\Parens1{\vrho-\lambda+2(\ell-1)}\gamma_{\ell-1}(\lambda)\\
      &+\Parens1{\vrho-\lambda+\ell-2}\Parens1{\vrho+\ell-2}\gamma_{\ell-2}(\lambda).
    \end{split}
  \end{equation}
  Conversely, define $\gamma_\ell(\lambda)$ by the above equation. Then the series $\Phi_\lambda $ converges to a solution of Equation \eqref{eq:phi_laplace}, absolutely on $A_\eps^+$ \fa $\eps>0$.
\end{Prop}

\begin{proof}
  In view of \thmref{Prop}{phi_laplace}, we may follow the standard procedure from Ref.~\cite{Hel84}*{Chapter IV, \S~5}. We obtain the recursion relation 
  \begin{equation}\label{eq:gamma-recursion}
    \ell(\ell-\lambda)\gamma_\ell(\lambda)=\sum_{n=1}^2\frac{nm_{n\alpha}}2\sum_{k\sge 1,\ell\sge nk}\Parens1{\vrho-\lambda+2(\ell-nk)}\gamma_{\ell-nk}(\lambda).
  \end{equation}
  On the right-hand side, the contribution of the three terms corresponding to $(n,k)=(1,1),(1,2),(2,1)$ is 
  \[
    \frac{m_\alpha}2\Parens1{\vrho-\lambda+2(\ell-1)}\gamma_{\ell-1}(\lambda)+\vrho(\vrho-\lambda+2(\ell-2))\gamma_{\ell-2}(\lambda).
  \]
  On applying the recursion relation to the remaining terms, we see that the remainder amounts to 
  \[
    (\ell-2)(\ell-2-\lambda)\gamma_{\ell-2}(\lambda).
  \]
  Thus, we obtain Equation \eqref{eq:coeff_Phi}.

  For the convergence, we follow the procedure of \cite{Hel84}*{Chapter IV, \S5, Lemma 5.3}: Since we have 
  \[
    \frac{\ell(\ell-\lambda)}{\ell^2}\longrightarrow 1,\quad\frac{\vrho-\lambda+2\ell}{\ell+1}\longrightarrow 2
  \]
  for $\ell\to\infty$, there are constants $0<c<\tfrac12$ and $2<C<\infty$ \scth 
  \[
    \Abs0{\vrho-\lambda+2\ell}\sle C(\ell+1),\quad\Abs0{\ell(\ell-\lambda)}\sge 2c\ell^2
  \]
  for all $\ell\in\nats$. For any $\eps>0$, there is some $\ell_0\in\nats$, $\ell_0\sge2$, \scth 
  \[
    \sum_{n=1}^2\frac{nm_{n\alpha}}4\Parens1{\coth\Parens1{\tfrac{nt}2}-1}\sle\frac{c\ell_0}C,\quad\forall t:t\sge\eps.
  \]
  
  Choose $K=K_\eps>0$ \scth 
  \[
    \Abs0{\gamma_\ell(\lambda)}\sle Ke^{\ell t},\quad\forall \ell,t:\ell\sle\ell_0,t\sge\eps.
  \]
  We claim that the `Gangolli estimate' $\Abs0{\gamma_\ell(\lambda)}\sle Ke^{\ell t}$ holds \fa $\ell\in\nats$ and $t\sge\eps$. Indeed, assume that this holds for all $\ell'<\ell_1$, where $\ell_1>\ell_0$. Then Equation \eqref{eq:gamma-recursion} and the assumption imply
  \begin{align*}
    \Abs0{\gamma_{\ell_1}(\lambda)}
    &\sle\frac{C(\ell_1+1)}{2c\ell_1^2}\sum_{n=1}^2\frac{n\Abs0{m_{n\alpha}}}2\sum_{k\sge 1,\ell\sge nk}\Abs0{\gamma_{\ell-nk}(\lambda)}\\
    &\sle\frac{Ke^{\ell t}C}{\ell_0c}\sum_{n=1}^2\frac{n\Abs0{m_{n\alpha}}}4\Parens1{\coth\Parens1{\tfrac{nt}2}-1}\sle Ke^{\ell t}    
  \end{align*}
  because $\ell_0(\ell_1+1)\sle 2\ell_1^2$ and $\coth(\tfrac t2)-1=2\sum_{k=1}^\infty e^{-kt}$. This proves our claim. 

  In particular, we have 
  \[
    \sum_{\ell=0}^\infty\Abs1{\gamma_\ell(\lambda)e^{-2\ell t}}\sle K\sum_{\ell=0}^\infty e^{-\ell t}\sle\frac K{1-e^{-\eps}}<\infty
  \]
  and the series converges absolutely on $A^+_\eps$. This justifies the term-by-term differentiation and implies that the limit $\Phi_\lambda$ is a solution of Equation \eqref{eq:phi_laplace}.
\end{proof}

\begin{Rem}
  When $m_{2\alpha}=0$ (\ie in the orthosymplectic case), the two-term recursion for $\gamma_\ell(\lambda)$ simplifies, and we easily obtain the following closed expression 
  \begin{equation}\label{eq:gamma-m2azero}
    \begin{split}
      \gamma_\ell(\lambda)&=\gamma _0(\lambda )\prod_{m=0}^{\ell-1}\frac{(m+\varrho )(m+\varrho -\lambda )}{(m+1)(m+1-\lambda )}\\
      &=\gamma _0(\lambda )c(-\lambda )(-1)^l\binom{-\varrho }{\ell}\frac{-\lambda}{(\ell-\lambda )c(\ell-\lambda )}.
    \end{split}
  \end{equation}
  Using this formula, the convergence of the Harish-Chandra series follows from a simple-minded application of the ratio test. 
\end{Rem}

From now on, will make the choice 
\begin{equation}\label{eq:gamma0-norm}
    \gamma_0(\lambda )\defi c(\lambda).
\end{equation}

Contrary to the classical situation of Riemannian symmetric spaces of non-compact type, the Harish-Chandra series may be finite. 

\begin{Cor}[sph_fin]
  If $G=\SOSp_{cs}(1,1+p|2q)$ and $m_\alpha=p-2q\leq 0$ is even, then the series $\Phi_\lambda$ is finite. More precisely, we have 
  \begin{align*}
    \Phi _\lambda (e^{th_0})=c_1e^{(\lambda-\vrho)t}\sum_{\ell=0}^{-\varrho }\binom{\lambda -\varrho }{\ell}\binom{-\lambda -\varrho }{-\varrho -\ell}e^{-2\ell t}.
  \end{align*}
  for some constant $c_1$. In particular, $\Phi_\lambda $ is well-defined on $A$. 
\end{Cor}

\begin{proof}
  In view of Equation \eqref{eq:gamma-m2azero}, we have $\gamma _\ell(\lambda )=0$ for $\ell>-\varrho $. Recall from \thmref{Prop}{c_osp_2} and \thmref{Cor}{c_osp_4} that
  \begin{equation}\label{eq:clambda-fin}
    c(\lambda )=c_0'\frac{\Gamma (\lambda )}{\Gamma (\lambda +\varrho )}=c_0'\prod_{k=1}^{-\varrho }(\lambda -k).
  \end{equation}
  Inserting $\gamma_0(\lambda)=c(\lambda)$ into Equation \eqref{eq:gamma-m2azero}, we obtain
  \begin{equation}\label{eq:gammaell-finite}
    \begin{split}
      \gamma_\ell(\lambda )
      &=c_0'(-1)^\ell\binom{-\varrho }{\ell}\prod_{k=\ell+1}^{-\varrho }(\lambda-k)\prod_{m=0}^{\ell-1}(\lambda -\varrho -m)\\
      &=c_0'(-\varrho )!\binom{-\lambda -\varrho }{-\varrho -\ell}\binom{\lambda -\varrho }{\ell}.    
    \end{split}    
  \end{equation}
  By the definition of $\Phi_\lambda$, this proves the claim.
\end{proof}

Combining the above results, we arrive at the following expression for the spherical superfunctions. 

\begin{Th}[sph_sum]
  Let $(G,K)$ be one of the symmetric pairs in Equation \eqref{eq:symmpairs-exgl11}. For $\Re\lambda>0$ and $\lambda\notin\tfrac12\ints$, we have
  \[
    \phi_\lambda|_{A^+}=\sum_{w\in W_0}\Phi_{w\lambda},
  \]
  where $\Phi_\lambda$ is defined by Equations \eqref{eq:def_Phi} and \eqref{eq:gamma0-norm}. Moreover, $W_0=\{\pm1\}$ unless $G=\SOSp_{cs}(1,1|2q)$, in which case $W_0=\{1\}$.
\end{Th}

\begin{proof}
  Firstly, Equation \eqref{eq:c_func}, \thmref{Th}{gk-fmla}, and \thmref{Prop}{Phi-conv} lead to 
  \begin{equation}\label{eq:phiPhi-asymp}
      \lim_{t\to\infty}\phi _{\lambda }(e^{th_0})e^{-(\lambda -\varrho )t}=c(\lambda )=\lim_{t\to\infty}\Phi _{\lambda }(e^{th_0})e^{-(\lambda -\varrho )t}.
  \end{equation}
  Therefore, in case $W_0$ is trivial, then \thmref{Prop}{spherical_osp_p0} and \thmref{Cor}{sph_fin} combine to prove $\phi_\lambda=\Phi_\lambda$.

  It remains to prove the assertion in case $W_0$ is non-trivial. We observe that by \thmref{Prop}{Phi-conv} again, the functions $\Phi_{\pm\lambda}$ are solutions of Equation \eqref{eq:phi_laplace}. This is a differential equation of order two and $\Phi_{\pm\lambda}$ are visibly linearly independent, so 
  \[
    \phi_\lambda=b_1\Phi_\lambda+b_2\Phi_{-\lambda },\quad b_1,b_2\in\cplxs.
  \]
  Let $w_0\in M_0'$ be a representative of the non-trivial Weyl group element. Using \thmref{Cor}{spherical_inv} and the $K\times K$-invariance of $\phi_\lambda $, we derive for any $a\in A^+$ that 
  \[
    \phi_\lambda(a)=\phi _{-\lambda }(a^{-1})=\phi _{-\lambda }(waw^{-1})=\phi _{-\lambda }(a).
  \] 
  This implies $b_1=b_2$. Using Equation \eqref{eq:phiPhi-asymp}, and 
  \[
    \lim_{t\to\infty}\Phi_{-\lambda}(e^{th_0})e^{-(\lambda-\vrho)t}=0,
  \] 
  which is derived easily from \thmref{Prop}{Phi-conv}, we conclude that $b_1=1$. 
\end{proof}

The result of \thmref{Th}{sph_sum} can be given a more explicit form in the cases where the Harish-Chandra series is finite. 

\begin{Cor}[sphfn-jacobi]
  Let $G=\SOSp_{cs}(1,1+p|2q)$ and $m_\alpha=p-2q\leq 0$ be even. Then 
  \begin{align*}
    \phi _\lambda (e^{th_0})=c_1e^{(\lambda -\varrho )t}P^{(-\lambda ,2\varrho -1)}_{-\varrho }(1-2e^{-2t}),
  \end{align*}
  where the  $P_n^{(a,b)}$ denote the Jacobi polynomials.  
\end{Cor}

\begin{proof}
  By \cite{emot-vol2}*{\S 10.8, (16)}, we have 
  \[
    P^{(a,b)}_n(x)=\binom{n+a}n{}_2F_1\Parens3{-n,1+a+b+n;a+1;\frac{x-1}2},
  \]
  where ${}_2F_1$ is the Gaussian hypergeometric function. Inserting the hypergeometric series, we obtain easily that 
  \[
    P^{(a,b)}_n(x)=\frac{(a+n)!}{n!(a+b+n)!}\sum_{\ell=0}^n\binom n\ell\frac{(a+b+n+\ell)!}{(a+\ell)!}\Parens3{\frac{x-1}2}^\ell,
  \]
  where we write $y!=\Gamma(y+1)$. Let 
  \[
    n\defi-\vrho=q-\frac p2,\quad a\defi-\lambda,\quad b\defi 2\vrho-1=p-2q-1.
  \]
  On applying Equation \eqref{eq:gammaell-finite} and \thmref{Cor}{sph_fin}, we obtain
  \[
    \Phi_\lambda(e^{th_0})=c_0'(-1)^{-\vrho}P^{(a,b)}_n(x)
  \]
  for $x=1-2e^{-2t}$, $c_0'$ denoting the constant from Equation \eqref{eq:clambda-fin}.
\end{proof}

\end{document}